\documentclass[twocolumn]{article}
\usepackage{fullpage,times}
\usepackage[final]{graphics}

\title{Numerical Study of Quantum Resonances in Chaotic Scattering}
\date{}

\author{Kevin K. Lin\thanks{Department of Mathematics, University of
California, Berkeley, CA 94720.  E-mail: {\bf kkylin@math.berkeley.edu}.
The author is supported by a fellowship from the Fannie and John Hertz
Foundation.}}

\begin{document}

  \newcommand{\Gerard}{G\'{e}rard}
  \newcommand{\Halpha}{\widehat{H}_{\alpha}}
  \newcommand{\Hhat}{\widehat{H}}
  \newcommand{\Hnalpha}{\widehat{H}_{N,\alpha}}
  \newcommand{\KEhatt}[1]{\widetilde{K}_E^{\of{{#1}}}}
  \newcommand{\KEhat}{\widetilde{K}_E}
  \newcommand{\Kint}[1]{K_{[{#1}]}}
  \newcommand{\PhiP}{\widetilde{\Phi}}
  \newcommand{\Poincare}{Poincar\'{e}}
  \newcommand{\Schrodinger}{Schr\"{o}dinger}
  \newcommand{\Sjostrand}{Sj\"{o}strand}
  \newcommand{\abs}[1]{\left|{#1}\right|}
  \newcommand{\imag}[1]{Im\left({#1}\right)}
  \newcommand{\interval}[1]{\left[{#1}\right]}
  \newcommand{\norm}[1]{\left|\left|{#1}\right|\right|}
  \newcommand{\of}[1]{\left({#1}\right)}
  \newcommand{\real}[1]{Re\left({#1}\right)}
  \newcommand{\set}[1]{\left\{{#1}\right\}}
  \newcommand{\vol}{\mbox{vol}}

  \newcommand{\PSbox}[3]{\includegraphics[0in,0in][{#2},{#3}]{#1}}

  \maketitle
  \begin{abstract}
  This paper presents numerical evidence that for quantum systems with
  chaotic classical dynamics, the number of scattering resonances near
  an energy $E$ scales like $\hbar^{-\frac{D\of{K_E}+1}{2}}$ as
  $\hbar\rightarrow{0}$.  Here, $K_E$ denotes the subset of the
  classical energy surface $\set{H=E}$ which stays bounded for all time
  under the flow generated by the Hamiltonian $H$ and $D\of{K_E}$
  denotes its fractal dimension.  Since the number of bound states in a
  quantum system with $n$ degrees of freedom scales like $\hbar^{-n}$,
  this suggests that the quantity $\frac{D\of{K_E}+1}{2}$ represents the
  effective number of degrees of freedom in scattering problems.
  \end{abstract}

\section{Introduction}
Quantum mechanics identifies the energies of stationary states in an
isolated physical system with the eigenvalues of its Hamiltonian
operator.  Because of this, eigenvalues play a central role in the study
of bound states, such as those describing the electronic structures of
atoms and molecules.\footnote{For examples, see {\cite{bethe}}.} When
the corresponding classical system allows escape to infinity, resonances
replace eigenvalues as fundamental quantities: The presence of a
resonance at $\lambda=E-i\gamma$, with $E$ real and $\gamma>0$, gives
rise to a dissipative metastable state with energy $E$ and decay rate
$\gamma$, as described in {\cite{z1}}.  Such states are essential in
scattering theory.\footnote{Systems which are not effectively isolated
but interact only weakly with their environment can also exhibit
resonant behavior.  For example, electronic states of an ``isolated''
hydrogen atom are eigenfunctions of a self-adjoint operator, but
coupling the electron to the radiation field turns those eigenstates
into metastable states with finite lifetimes.  This paper does not deal
with dissipative systems and is only concerned with scattering.}

An important property of eigenvalues is that one can count them using
only the classical Hamiltonian function $H(x,p)=\frac{1}{2}\norm{p}^2 +
V(x)$ and Planck's constant $\hbar$: For fixed energies $E_0<E_1$, the
number $N_{eig}(E_0,E_1,\hbar)$ of eigenvalues in $\interval{E_0,E_1}$
is
\begin{equation}
  N_{eig}(E_0,E_1,\hbar) \approx \frac{\vol\of{\set{E_0 \leq H \leq
                                 E_1}}} {(2\pi\hbar)^n},
\end{equation}
where $n$ denotes the number of degrees of freedom and $\vol\of{\cdot}$
phase space volume.  This result, known as the {\em Weyl law}, expresses
the density of quantum states using the classical Hamiltonian
function.\footnote{For a beautiful exposition of early work on this and
related themes, see {\cite{kac}}.  For recent work in the semiclassical
context, see {\cite{DS}}.} No direct generalization to resonances is
currently known.

In this paper, numerical evidence for a Weyl-like power law is presented
for resonances in a two-dimensional model with three
symmetrically-placed gaussian potentials.  A conjecture, based on the
work of {\Sjostrand\ }{\cite{Sj}} and Zworski {\cite{z2}}, states that
the number of resonances $\lambda=E-i\gamma$ with $E_0<E<E_1$ and
$0<\gamma<\hbar$ asymptotically lies between
$C_1\hbar^{-\frac{D\of{K_{E_1}}+1}{2}}$ and
$C_0\hbar^{-\frac{D\of{K_{E_0}}+1}{2}}$ as $\hbar\rightarrow{0}$, where
\begin{equation}
  \begin{array}{ccl}
  D\of{\cdot} &=& \mbox{dimension (see below)},\\
  K_E &=& K \cap \set{H=E},\\
  K &=& \set{(x,p):\sup_{t}\norm{\Phi_t(x,p)}<\infty},\\
  \Phi_t &=& \mbox{flow generated by $H$}.\\
  \end{array}
\end{equation}
If $D\of{K_E}$ depends continuously on $E$ and $\abs{E_1-E_0}$ is
sufficiently small, then $D\of{K_{E_1}}\approx{D\of{K_{E_0}}}$ and the
number of resonances in such a region is comparable to
$\hbar^{-\frac{D\of{K_E}+1}{2}}$ for any $E\in\interval{E_0,E_1}$.

The sets $K$ and $K_E$ are {\em trapped sets} and consist of initial
conditions which generate trajectories that stay bounded forever.  In
systems where $\set{H\leq{E}}$ is bounded for all $E$, the conjecture
reduces to the Weyl asymptotic $\hbar^{-n}$.

The notion of dimension requires some comment: The ``triple gaussian''
model considered here has very few trapped trajectories, and $K$ and
$K_E$ (for any energy $E$) have vanishing Lebesgue measures.  Thus,
$D(K)$ is strictly less than $2n=4$ and $D\of{K_E}<2n-1=3$.  In fact,
the sets $K$ and $K_E$ are fractal, as are trapped sets in many other
chaotic scattering problems.  Also, in this paper, the term ``chaotic''
always means {\em hyperbolic}; see {\Sjostrand\ }\cite{Sj} or Gaspard
{\cite{gaspard}} for definitions.

This paper is organized as follows: First, the model system is defined.
This is followed by mathematical background information, as well as a
heuristic argument for the conjecture.  Then, numerical methods for
computing resonances and fractal dimensions are developed, and numerical
results are presented and compared with known theoretical predictions.

\paragraph{Notation.}
In this paper, $H$ denotes the Hamiltonian function
$\frac{1}{2}\norm{p}^2 + V(x)$ and $\Hhat$ the corresponding Hamiltonian
operator $-\frac{\hbar^2}{2}\Delta + V(x)$, where
$\Delta=\sum_{k=1}^{n}{\frac{\partial^2}{\partial{x^2_k}}}$ is the usual
Laplacian and $V$ acts by multiplication.

\section{Triple Gaussian Model}
The model system has $n=2$ degrees of freedom; its phase space is $R^4$,
whose points are denoted by $(x,y,p_x,p_y)$.

First, it is convenient to define
\begin{equation}
  G_{x_0}^\sigma(x) = \exp\of{-\frac{(x-x_0)^2}{2\sigma^2}}.
\end{equation}
Similarly, put
\begin{equation}
  \begin{array}{lcl}
  G_{(x_0,y_0)}^\sigma(x,y) &=&
  \of{G_{x_0}^\sigma\otimes{G_{y_0}^\sigma}}(x,y)\\
  &=& G_{x_0}^\sigma(x)\cdot{G_{y_0}^\sigma(y)}\\
  \end{array}
\end{equation}
in two dimensions.

Now, define $H$ by
\begin{equation}
  \label{eqn:classical-hamiltonian}
  H(x,y,p_x,p_y) = \frac{1}{2}\of{p_x^2 + p_y^2} + V_3(x,y)
\end{equation}
where the potential $V_m$ is given by
\begin{equation}
  \label{eqn:potdef}
  \begin{array}{ccl}
  V_m &=& \sum_{k=1}^m{G_{c(k,m)}^\sigma},\\
  c(k,m) &=& (R\cos\of{\frac{2\pi k}{m}},R\sin\of{\frac{2\pi k}{m}}).\\
  \end{array}
\end{equation}
That is, it consists of $m$ gaussian ``bumps'' placed at the vertices of
a regular $m$-gon centered at the origin, at a distance $R>0$ from the
origin.  This paper focuses on the case $m=3$ because it is the simplest
case that exhibits nontrivial dynamics in two dimensions.  However, the
case $m=2$ is also relevant because it is well-understood: See Miller
{\cite{miller0}} for early heuristic results and {\Gerard\ }and
{\Sjostrand\ }{\cite{GS}} for a rigorous treatment.  Thus, double
gaussian scattering serves as a useful test case for the techniques
described here.

\begin{figure}
  \begin{center}
  \PSbox{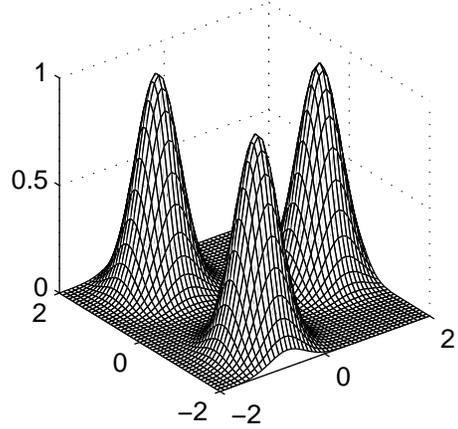}{2.5in}{2.5in}
  \end{center}
  \caption{Triple gaussian potential}
  \label{fig:bumps}
\end{figure}

The quantized Hamiltonian $\Hhat$ is similarly defined:
\begin{equation}
  \Hhat=-\frac{\hbar^2}{2}\Delta + V_3.
\end{equation}
See Figure {\ref{fig:bumps}}.

\section{Background}
This section provides a general discussion of resonances and motivates
the conjecture in the context of the triple gaussian model.  However,
the notation reflects the fact that most of the definitions and
arguments here carry over to more general systems with $n$ degrees of
freedom.  The reader should keep in mind that $n=2$ for the triple
gaussian model.

There exists an extensive literature on resonances and semiclassical
asymptotics in other settings.  For example, see
{\cite{GR1,GR2,GR3,wirzba}} for detailed studies of the classical and
quantum mechanics of hard disc scattering.

\subsection{Resonances}
Resonances can be defined mathematically as follows: Set
$R(z)=(\Hhat-zI)^{-1}$ for real $z$, where $I$ is the identity operator.
This one-parameter family of operators $R(z)$ is the {\em resolvent} and
is meromorphic with suitable modifications of its domain and range.  The
poles of its continuation into the complex plane are, by definition, the
{\em resonances} of $\Hhat$.\footnote{For more details and some
references, see {\cite{z1}}.}

Less abstractly, resonances are generalized eigenvalues of $\Hhat$.
Thus, we should solve the time-independent {\Schrodinger\ }equation
\begin{equation}
  \label{eqn:time-indep-schrodinger}
  \Hhat\psi = \lambda\psi
\end{equation}
to obtain the resonance $\lambda$ and its generalized eigenfunction
$\psi$.  In bound state computations, one approximates $\psi$ as a
finite linear combination of basis functions and solves a
finite-dimensional version of the equation above.  To carry out similar
calculations for resonances, it is necessary that $\psi$ lie in a
function space which facilitates such approximations, for example $L^2$.

Let $\psi$ and $\lambda$ solve (\ref{eqn:time-indep-schrodinger}).  Then
$e^{-\frac{i} {\hbar}\lambda{t}}\cdot{\psi}$ solves the time-dependent
{\Schrodinger\ }equation
\begin{equation}
  \label{eqn:time-dep-schrodinger}
  i\hbar \frac{\partial\psi}{\partial{t}} = \Hhat{\psi}.
\end{equation}
It follows that $\imag{\lambda}$ must be negative because metastable
states decay in time.  Now suppose, for simplicity, that
$n=1$.\footnote{The analysis in higher dimensions requires some care,
but the essential result is the same.} Then solutions of
(\ref{eqn:time-indep-schrodinger}) with energy $E$ behave like
$e^{-\frac{i}{\hbar}\sqrt{2E}x}$ for large $x>0$.  Substituting
$\lambda=E-i\gamma$ for $E$ yields $e^{-\frac{i}
{\hbar}\sqrt{2\lambda}x}$, which grows exponentially because
$\imag{\sqrt{E-i\gamma}}<0$.  Thus, finite rank approximations of
$\Hhat$ cannot capture such generalized eigenfunctions.  However, if we
make the formal substitution $x\mapsto{xe^{i\alpha}}$, then the wave
function becomes $\exp\of{-\frac{i}
{\hbar}\sqrt{2\lambda}\cdot{e^{i\alpha}}\cdot{x}}$.  Choosing
$\alpha>\frac{1}{2}\tan^{-1}\of{\gamma/E}$ forces $\psi$ to decay
exponentially.

This procedure, called {\em complex scaling}, transforms the Hamiltonian
operator $\Hhat$ into the {\em scaled operator} $\Halpha$.  It also maps
metastable states $\psi$ with decay rate $\gamma<E\tan\of{2\alpha}$ to
genuine $L^2$ eigenfunctions $\psi_\alpha$ of $\Halpha$.  The
corresponding resonance $\lambda$ becomes a genuine eigenvalue:
$\Halpha\psi_\alpha=\lambda\psi_\alpha$.  Furthermore, resonances of
$\Hhat$ will be invariant under small perturbations in $\alpha$, whereas
other eigenvalues of $\Halpha$ will not.  The condition
$\alpha>\frac{1}{2}\tan^{-1}\of{\gamma/E}$ implies that, for small
$\gamma$ and fixed $E$, the method will capture a resonance
$\lambda=E-i\gamma$ if and only if $\gamma<2E\alpha+O\of{\alpha^2}$.  We
can perform complex scaling in higher dimensions by substituting
$r\mapsto{re^{i\alpha}}$ in polar coordinates.

\begin{figure}
  \begin{center}
  \PSbox{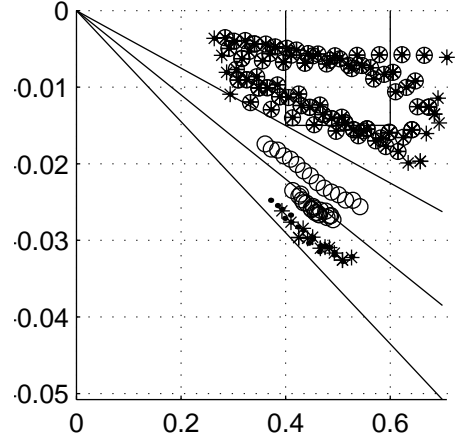}{2.5in}{2.5in}
  \end{center}

  \caption{Illustration of complex scaling: The three lines indicate the
  location of the rotated continuous spectrum for different values of
  $\alpha$, while the box at the top of the figure is the region in
  which resonances are counted.  Eigenvalues which belong to different
  values of $\alpha$ are marked with different styles of points.  As
  explained later, only eigenvalues near the region of interest are
  computed.  This results in a seemingly empty plot.}

  \label{fig:alpha}
\end{figure}

In algorithmic terms, this means we can compute eigenvalues of $\Halpha$
for a few different values of $\alpha$ and look for invariant values, as
demonstrated in Figure {\ref{fig:alpha}}.  In addition to its accuracy
and flexibility, this is one of the advantages of complex scaling: The
invariance of resonances under perturbations in $\alpha$ provides an
easy way to check the accuracy of calculations, mitigating some of the
uncertainties inherent in computational work.\footnote{For a different
approach to computing resonances, see {\cite{WMS}} and the references
there.} Note that the scaled operator $\Halpha$ is no longer
self-adjoint, which results in non-hermitian finite-rank approximations
and complex eigenvalues.

This method, first introduced for theoretical purposes by Aguilar and
Combes {\cite{AC}} and Balslev and Combes {\cite{BC}}, was further
developed by B.\ Simon in {\cite{simon}}.  It has since become one of
the main tools for computing resonances in physical chemistry
{\cite{miller1,miller2,miller3,McC}}.  For recent mathematical progress,
see {\cite{LM,Sj,SZ}} and references therein.

For reference, the scaled triple-gaussian operator $\Halpha$ is
\begin{equation}
\label{eqn:scaled-hamiltonian}
  \Halpha = -e^{-2i\alpha}\cdot\frac{\hbar^2}{2}\Delta + V_{3,\alpha},
\end{equation}
where
\begin{equation}
\label{eqn:scaled-potential}
  \begin{array}{ccl}
  V_{m,\alpha} &=&
  \sum_{k=1}^m{G_{c_\alpha(k,m)}^{\sigma_\alpha}},\\
  \sigma_\alpha &=& e^{-i\alpha}\cdot\sigma,\\
  c_\alpha(k,m) &=& e^{-i\alpha}c(k,m).\\
  \end{array}
\end{equation}
Note that these expressions only make sense because $G_{x_0}^\sigma(x)$
is analytic in $x$, $x_0$, and $\sigma$.

\subsection{Fractal Dimension}
Recall that the {\em Minkowski dimension} of a given set $U\subset{R^m}$
is
\begin{equation}
  \label{eqn:minkowski-def0}
  D=\inf\set{d:\limsup_{\epsilon\rightarrow{0}}\of{
         {\epsilon^{d-m} \cdot \vol\of{U_\epsilon}}} < \infty},
\end{equation}
where $U_\epsilon=\set{y\in{R^m}:\mbox{dist}(y,U)<\epsilon}$.  A simple
calculation yields
\begin{equation}
  \label{eqn:minkowski-def1}
  D(U)=\lim_{\epsilon\rightarrow{0}}
       \frac{\log\of{\vol\of{U_\epsilon}/\epsilon^m}}
       {\log\of{1/\epsilon}}
\end{equation}
if the limit exists.

Texts on the theory of dimensions typically begin with the Hausdorff
dimension because it has many desirable properties.  In contrast, the
Minkowski dimension can be somewhat awkward: For example, a countable
union of zero-dimensional sets (points) can have positive Minkowski
dimension.  But, the Minkowski dimension is sometimes easier to
manipulate and almost always easier to compute.  It also arises in the
heuristic argument given below.

For a detailed treatment of different definitions of dimension and their
applications in the study of dynamical systems, see
{\cite{falconer,pesin}}.

\subsection{Generalizing the Weyl Law}
The formula
\begin{equation}
  \frac{\vol\of{\set{E_0 \leq H \leq E_1}}}{(2\pi\hbar)^n}
\end{equation}
makes no sense in scattering problems because the volume on the right
hand side is infinite for most choices of $E_0$ and $E_1$, and this
seems to mean that there is no generalization of the Weyl law in the
setting of scattering theory.  However, the following heuristic argument
suggests otherwise:

As mentioned before, a metastable state corresponding to a resonance
$\lambda=E-i\gamma$ has a time-dependent factor of the form
$e^{-\frac{i}{\hbar}\lambda{t}} =
e^{-\frac{i}{\hbar}Et}\cdot{e^{-\frac{\gamma}{\hbar}{t}}}$.  A wave
packet whose dynamics is dominated by $\lambda$ (and other resonances
near it) would therefore exhibit temporal oscillations of frequency
$O(E/\hbar)$ and lifetime $O(\hbar/\gamma)$.  Heuristically, then, the
number of times the particle ``bounces'' in the ``trapping
region''\footnote{For our triple gaussian system, that would be the
triangular region bounded by the gaussian bumps.} before escaping should
be comparable to $\frac{E} {\hbar}\cdot\frac{\hbar} {\gamma} =\frac{E}
{\gamma}$.

In the semiclassical limit, the dynamics of the wave packet should be
well-approximated by a classical trajectory.  Let $T(x,y,p_x,p_y)$
denote the time for the particle to escape the system starting at
position $(x,y)$ with momentum $(p_x,p_y)$.  The diameter of the
trapping region is $O(R)$, and typical velocities in the energy surface
$\set{H=E}$ are $O(\sqrt{E})$ (mass set to unity), so the number of
times a classical particle bounces before escaping should be
$O(T\sqrt{E}/R)$.  This suggests that, in the limit
$\hbar\rightarrow{0}$, $T\sqrt{E}/R\sim{E/\gamma}$ and consequently
\begin{equation}
  \label{eqn:T-gamma}
  T \sim \frac{R\sqrt{E}}{\gamma}.
\end{equation}
Fix $\gamma_0>0$, and consider
\begin{equation}
  \label{eqn:Nres-def}
   N_{res}=\#\set{E-i\gamma:E_0 \leq E \leq E_1, \gamma\leq\gamma_0}
\end{equation}
for fixed energies $E_0$ and $E_1$: Equation (\ref{eqn:T-gamma}) implies
that $T\geq R\sqrt{E_0}/\gamma_0$, so by analogy with the Weyl law,
\begin{equation}
  \frac{\vol\of{\set{E_0 \leq H \leq E_1,
       T \geq \frac{R\sqrt{E_0}}{\gamma_0}}}} {(2\pi\hbar)^n}
\end{equation}
follows as an approximation for the number of quantum states with the
specified energies and decay rates.

Now, the function $1/T$ is nonnegative for all $(x,y,p_x,p_y)$ and
vanishes on $\Kint{E_0,E_1} = K\cap\set{E_0\leq{H}\leq{E_1}}$.  Assuming
that $1/T$ is sufficiently regular,\footnote{\label{fn:quad} In fact,
this is numerically self-consistent: Assume that $1/T$ vanishes to order
$\nu$ (with $\nu$ not necessarily equal to $2$) on $K$, and {\em assume}
the conjecture.  Then the number of resonances would scale like
$\hbar^{(2n-D(K))/\nu}$, from which one can solve for $\nu$.  With the
numerical data we have, this indeed turns out to be $2$ (but with
significant fluctuations).

Also, if $1/T$ does not vanish quadratically {\em everywhere} on $K$,
variations in its regularity may affect the correspondence between
classical trapping and the distribution of resonances.} this suggests
\begin{equation}
   1/T(x,y,p_x,p_y) \sim d_{\Kint{E_0,E_1}}(x,y,p_x,p_y)^2,
\end{equation}
where $d_{\Kint{E_0,E_1}}$ denotes distance to $\Kint{E_0,E_1}$.  It
follows that $N_{res}$ should scale like
\begin{equation}
  \frac{\vol\of{\set{E_0 \leq H \leq E_1, d_{\Kint{E_0,E_1}} \leq
  \gamma_0^\frac{1}{2}}}} {\hbar^n}.
\end{equation}
For small $\gamma_0$, this becomes
\begin{equation}
  C(R,E_0,E_1)\cdot\hbar^{-n} \cdot
  \gamma_0^{\frac{2n-D\of{\Kint{E_0,E_1}}}{2}}
\end{equation}
for some constant $C$, by (\ref{eqn:minkowski-def0}).  Choosing
$\gamma_0=\hbar$ and assuming that $D\of{K_E}$ decreases monotonically
with increasing $E$ (as is the case in Figure {\ref{fig:dimconv}}), we
obtain
\begin{equation}
  C_1\hbar^{-\frac{D\of{K_{E_1}}+1}{2}} \leq N_{res} \leq
  C_0\hbar^{-\frac{D\of{K_{E_0}}+1}{2}}.
\end{equation}
If $\abs{E_1-E_0}$ is sufficiently small, then
$D\of{\Kint{E_0,E_1}}\approx{D\of{K_E}+1}$ for $E\in\interval{E_0,E_1}$,
and
\begin{equation}
  \label{eqn:scaling-law}
  N_{res} \sim \hbar^{-\frac{D\of{K_E}+1}{2}}.
\end{equation}
In {\cite{Sj}}, {\Sjostrand\ }proved the following rigorous upper bound:
For $\gamma_0>0$ satisfying $C\hbar<\gamma_0<1/C$,
\begin{equation}
  N_{res} =
  O\of{C_\delta\hbar^{-n}\gamma_0^{\frac{2n-D\of{\Kint{E_0,E_1}} +
  \delta}{2}}}
\end{equation}
holds for all $\delta>0$.  When the trapped set is of {\em pure
dimension}, that is when the infimum in Equation
(\ref{eqn:minkowski-def0}) is achieved, one can take $\delta=0$.
Setting $\gamma_0=\hbar$ gives an upper bound of the form
(\ref{eqn:scaling-law}).

In his proof, {\Sjostrand\ }used the semiclassical argument above with
escape functions and the Weyl inequality for singular values.  Zworski
continued this work in {\cite{z2}}, where he proved a similar result for
scattering on convex co-compact hyperbolic surfaces with no cusps.  His
work was motivated by the availability of a large class of examples with
hyperbolic flows, easily computable dimensions, and the hope that the
Selberg trace formula could help obtain lower bounds.  But, these hopes
remain unfulfilled so far {\cite{GZ}}, and that partly motivates this
work.

\section{Computing Resonances}
Complex scaling reduces the problem of calculating resonances to one of
computing eigenvalues.  What remains is to approximate the operator
$\Halpha$ by a rank $N$ operator $\Hnalpha$ and to develop appropriate
numerical methods.  For comparison, see
{\cite{miller1,miller2,miller3,McC}} for applications of complex scaling
to problems in physical chemistry.

\subsection{Choice of Scaling Angle.}
\label{sec:alpha-choice}
One important consideration in resonance computation is the choice of
the scaling angle $\alpha$: Since we are interested in counting
resonances in a box $\interval{E_0,E_1}-i\interval{0,\hbar}$, it is
necessary to choose $\alpha\geq{\tan^{-1}\of{\frac{\hbar}{E_0}}}$ so
that the continuous spectrum of $\Halpha$ is shifted out of the box
$\interval{E_0,E_1}-i\interval{0,\hbar}$ (see Figure {\ref{fig:alpha}}).

In fact, the resonance calculation uses
\begin{equation}
  \begin{array}{lcl}
    \alpha &=& \tan^{-1}\of{\frac{\hbar}{E_0}}\\
           &=& \frac{\hbar}{E_0} + O\of{\hbar^2}.\\
  \end{array}
\end{equation}
This choice of $\alpha$ helps avoid the {\em pseudospectrum}
{\cite{trefethen,z3}}:

Let $A$ be an $N\times{N}$ matrix, and let $R(z)$ be the resolvent
$(A-zI)^{-1}$.  It is well known that when $A$ is {\em normal}, that is
when $A$ commutes with its adjoint $A^*$, the spectral theorem applies
and the inequality
\begin{equation}
  \begin{array}{lcl}
  \norm{R(z)} &=& \norm{(A-zI)^{-1}}\\
           &\leq& \mbox{dist}(z,\sigma(A))^{-1}\\
  \end{array}
\end{equation}
holds ($\sigma(A)$ denotes the spectrum of $A$).  When $A$ is not
normal, no such inequality holds and $\norm{R(z)}$ can become very large
for $z$ far from $\sigma(A)$.  This leads one to define {\em
$\epsilon$-pseudospectrum}:
\begin{equation}
  \Lambda_\epsilon(A) = \set{z:\norm{R(z)} \geq 1/\epsilon}.
\end{equation}
Using the fact that $A$ is a matrix, one can show that
$\Lambda_\epsilon(A)$ is equal to the set
\begin{equation}
  \set{z:\exists{A'} \mbox{ such that } z\in\sigma(A+A'),
  \norm{A'}\leq\epsilon}.
\end{equation}
That is, the $\epsilon$-pseudospectrum of $A$ consists of those complex
numbers $z$ which are eigenvalues of an $\epsilon$-perturbation of $A$.

\begin{figure}
  \begin{center}
  \PSbox{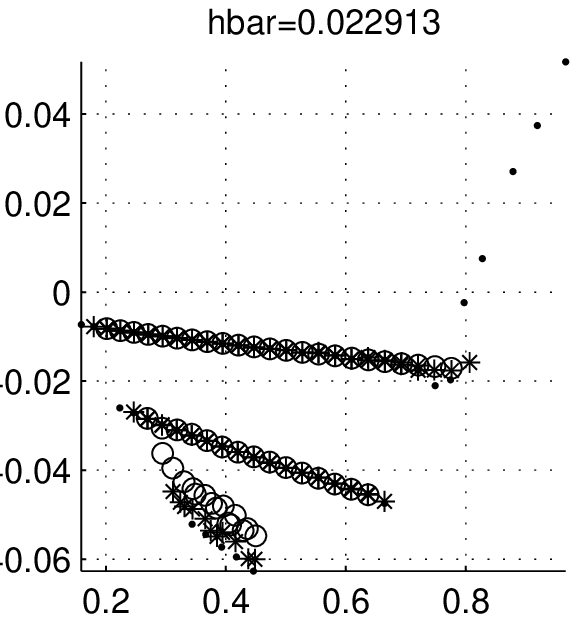}{2.5in}{2.5in}
  \end{center}
  \caption{Resonances for two-bump scattering with $\hbar=0.022913$.}
  \label{fig:twobumps3}
\end{figure}

The idea of pseudospectrum can be extended to general linear operators.
In {\cite{trefethen}}, it is emphasized that for non-normal operators,
the pseudospectrum can create ``false eigenvalues'' which make the
accurate numerical computation of eigenvalues difficult.  In
{\cite{z3}}, this phenomenon is explained using semiclassical
asymptotics.  Roughly speaking, the pseudospectrum of the scaled
operator $\Halpha$ is given by the closure of
\begin{equation}
  \set{z:z=H_\alpha\of{x,y,p_x,p_y}}
\end{equation}
of its {\em symbol} $H_\alpha$, which is the scaled Hamiltonian function
\begin{equation}
  H_\alpha(x,y,p_x,p_y) = e^{-2i\alpha}\cdot\frac{p_x^2+p_y^2}{2} +
  V_{3,\alpha}
\end{equation}
in this case.  Choosing $\alpha$ to be comparable to $\hbar$ ensures
that the imaginary part of $H_\alpha$ is also comparable to $\hbar$,
which keeps the pseudospectrum away from the counting box
$\interval{E_0,E_1}-i\interval{0,\hbar}$; a larger $\alpha$ would
contribute a larger $\alpha^2$ term to the imaginary part of $H_\alpha$
and enlarge the pseudospectrum.  As one can see in Figures
{\ref{fig:twobumps0}} - {\ref{fig:twobumps3}}, the invariance of
resonances under perturbations in $\alpha$ also helps filter out
pseudospectral effects.

This consideration also points out the {\em necessity} of the choice
$\alpha=\tan^{-1}\of{\frac{\hbar}{E_0}}$: To avoid pseudospectral
effects, $\alpha$ must be $O\of{\hbar}$.  On the other hand, if
$\alpha=o\of{\hbar}$, then finite rank approximations may fail to
capture resonances in the region of interest.

\subsection{Eigenvalue Computation}
Suppose that we have constructed $\Hnalpha$.  In the case of
eigenvalues, the Weyl law states that $N_{eig}=O\of{\hbar^{-2}}$ as
$\hbar\rightarrow{0}$, since our system has $n=2$ degrees of freedom.
Thus, in order to capture a sufficient number of eigenvalues, the rank
$N$ of the matrix approximation must scale like $\hbar^{-2}$.  In the
absence of more detailed information on the density of resonances, the
resonance computation requires a similar assumption to ensure sufficient
numerical resolution.

Thus, for moderately small $\hbar$, the matrix has $N^2\sim\hbar^{-4}$
entries, which rapidly becomes prohibitive on most computers available
today.  Furthermore, even if one does not store the entire matrix,
numerical packages like LAPACK {\cite{lapack}} require $O(N^2)$
auxiliary storage, again making practical calculations impossible.

Instead of solving the eigenvalue problem $\Hnalpha{v}=\lambda{v}$
directly, one solves the equivalent eigenvalue problem
\begin{equation}
  \of{\Hnalpha - \lambda_0}^{-1}v=\lambda'v.
\end{equation}
Efficient implementations of the Arnoldi algorithm {\cite{arpack}} can
solve for the largest few eigenvalues $\lambda'$ of
$\of{\Hnalpha-\lambda_0}^{-1}$.  But $\lambda=\lambda_0+1/\lambda'$, so
this method allows one to compute a subset of the spectrum of $\Hnalpha$
near a given $\lambda_0$.

Such algorithms require a method for applying the matrix
$\of{\Hnalpha-\lambda_0}^{-1}$ to a given vector $v$ at each iteration
step.  In the resonance computation, this is done by solving
$(\Hnalpha-\lambda_0)w=v$ for $w$ by applying conjugate gradient to the
normal equations (see {\cite{templates}}).\footnote{That is, instead of
solving $Aw=v$, one solves $A^*Aw=A^*v$.  This is necessary because
$\Hnalpha$ is non-hermitian, and conjugate gradient only works for
positive definite matrices.  This is not the best numerical method for
non-hermitian problems, but it is easy to implement and suffices in this
case.} The resonance program, therefore, consists of two nested
iterative methods: An outer Arnoldi loop and an inner iterative linear
solver for $(\Hnalpha-\lambda_0)w=v$.  This computation uses
ARPACK\footnote{See {\cite{arpack}} for details on the package, as well
as an overview of Krylov subspace methods.}, which provides a flexible
and efficient implementation of the Arnoldi method.

To compute resonances near a given energy $E$, the program uses
$\lambda_0=E+ia$, $a>0$, instead of $\lambda_0=E$: This helps control
the condition number of $\Hnalpha - \lambda_0$ and gives better error
estimates and convergence criteria.\footnote{Most of the error in
solving the matrix equation $(\Hnalpha-\lambda_0)w=v$ concentrates on
eigenspaces of $(\Hnalpha - \lambda_0)^{-1}$ with large eigenvalues.
These are precisely the desired eigenvalues, so in principle one can
tolerate inaccurate solutions.  However, the calculation requires
convergence criteria and error estimates for the linear solver, and
using $a>0$, say $a=1$, turns out to ensure a relative error of about
$10^{-6}$ after about 17-20 iterations of the conjugate gradient solver.
Since we only wanted to count eigenvalues, a more accurate (and
expensive) computation of resonances was not necessary.}

\subsection{Matrix Representations}
\subsubsection{Choice of Basis}
While one can discretize the differential operator $\Halpha$ via finite
differences, in practice it is better to represent the operator using a
basis for a subspace of the Hilbert space $L^2$: This should better
represent the properties of wave functions near infinity and obtain
smaller (but more dense) matrices.

\begin{figure}
  \begin{center}
  \PSbox{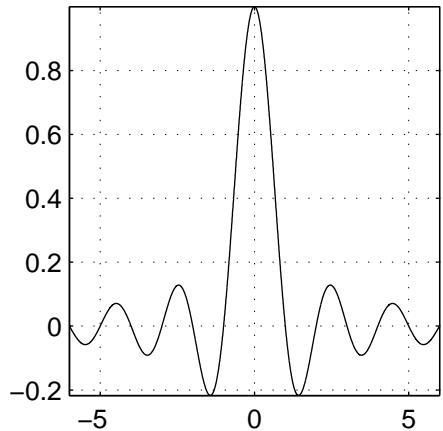}{2.5in}{2.5in}
  \end{center}
  \caption{A sinc function with $m=0$, $\Delta{x}=1$.}
  \label{fig:sinc}
\end{figure}

Common basis choices in the chemical literature include so-called
``phase space gaussian'' {\cite{DH}} and ``distributed gaussian'' bases
{\cite{HL}}.  These bases are not orthogonal with respect to the usual
$L^2$ inner product, so one must explicitly orthonormalize the basis
before computing the matrix representation of $\Halpha$.  In addition to
the computational cost, this also requires storing the entire matrix and
severely limits the size of the problem one can solve.  Instead, this
computation uses a {\em discrete-variable representation} (DVR) basis
{\cite{LHL}}:

Consider, for the moment, the one dimensional problem of finding a basis
for a ``good'' subspace of $L^2(R)$.  Fix a constant $\Delta{x}>0$, and
for each integer $m$, define
\begin{equation}
  \phi_{m,\Delta{x}}(x) =
  \sqrt{\Delta{x}}\cdot\frac{\sin\of{\frac{\pi}{\Delta{x}}(x -
  m\Delta{x})}}{\pi(x - m\Delta{x})}.
\end{equation}
(This is known as a ``sinc'' function in engineering literature
{\cite{siebert}}.  See Figure {\ref{fig:sinc}}.)  The Fourier transform
of $\phi_{m,\Delta{x}}$ is
\begin{equation}
  \widehat{\phi}_{m,\Delta{x}}(\omega) = \left\{
    \begin{array}{ll}
      e^{-im\Delta{x}}\cdot\sqrt{\Delta{x}},&\abs{\omega}\leq\pi/\Delta{x}\\
      0,&\abs{\omega}>\pi/\Delta{x}\\
    \end{array}\right.
\end{equation}
One can easily verify that $\set{\phi_{m,\Delta{x}}}$ forms an
orthonormal basis for the closed subspace of $L^2$ functions whose
Fourier transforms are supported in
$\interval{-\pi/\Delta{x},\pi/\Delta{x}}$.

\begin{figure}
  \begin{center}
  \PSbox{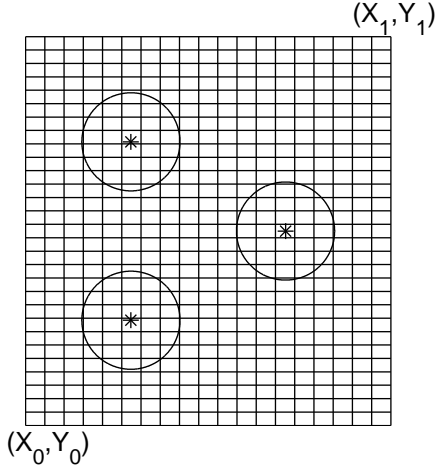}{2.5in}{2.5in}
  \end{center}

  \caption{Illustration of resonance program parameters in configuration
  space: The lower-left corner of the mesh is $(X_0,Y_0)$, while the
  upper right corner is $(X_1,Y_1)$.  The mesh contains $N_x\times{N_y}$
  grid points, and a basis function $\phi_{mn}$ is placed at each grid
  point.  Stars mark the centers of the potentials, the circles have
  radius $2\sigma$ (with $\sigma=1/3$), and $R$ is set to $1.4$.
  Parameters for the classical computation are depicted in Figure
  {\ref{fig:dim-prog-params}}.}

  \label{fig:res-prog-params}
\end{figure}

To find a basis for corresponding space of band-limited functions in
$L^2(R^2)$, simply form the tensor products
\begin{equation}
  \phi_{mn}(x_0,y_0) = \phi_{m,\Delta_x}(x)\phi_{n,\Delta_y}(y).
\end{equation}
The basis has a natural one-to-one correspondence with points
$(m\Delta_1+X_0,n\Delta_2+Y_0)$ on a regular lattice of grid points in a
box $\interval{X_0,X_1}\times\interval{Y_0,Y_1}$ covering the spatial
region of interest.  (See Figure {\ref{fig:res-prog-params}}.)  Using
this basis, it is easy to compute matrix elements for $\Halpha$.

\subsubsection{Tensor Product Structure}
An additional improvement comes from the separability of the
Hamiltonian: Each term in the scaled Hamiltonian $\Halpha$ splits into a
tensor product:
\begin{eqnarray}
  \frac{\partial^2}{\partial{x^2}} &=& \frac{d^2}{dx^2}\otimes{I_y}\\
  \frac{\partial^2}{\partial{y^2}} &=& {I_x}\otimes\frac{d^2}{dy^2}\\
  G_{(x_0,y_0)}^\sigma &=&
  G_{x_0}^\sigma\otimes{G_{y_0}^\sigma},
\end{eqnarray}
where $I_x$ and $I_y$ denote identity operators on copies of $L^2(R)$.
Since the basis $\set{\phi_{mn}}$ consists of tensor products of one
dimensional bases, $\Hnalpha$ is also a short sum of tensor products.
Thus, if we let $N_x$ denote the number of grid points in the $x$
direction and let $N_y$ denote the number of grid points in the $y$
direction, then $N=N_x\cdot{N_y}$ and $\Hnalpha$ is a sum of five
matrices of the form $A_x\otimes{A_y}$, where $A_x$ is $N_x\times{N_x}$
and $A_y$ is $N_y\times{N_y}$.

Such tensor products of matrices can be applied to arbitrary vectors
efficiently using the outer product representation.\footnote{The tensor
product of two column vectors $v$ and $w$ can be represented as
$v\cdot{w^T}$.  We then have
$(A\otimes{B})\cdot(v\otimes{w})=(Av)\cdot(Bw)^T$, which extends by
linearity to $(A\otimes{B})\cdot{u}=A\cdot{u}\cdot{B^T}$.} Since the
rank of $\Hnalpha$ is $N=N_x\cdot{N_y}$ and $N_x\approx{N_y}$ in these
computations, we can store the tensor factors of the matrix $\Hnalpha$
using $O(N)$ storage instead of $O(N^2)$, and apply $\Hnalpha$ to a
vector in $O\of{N^{3/2}}$ time instead of $O\of{N^3}$.  The resulting
matrix is not sparse, as one can see from the matrix elements for the
Laplacian below.

Note that this basis fails to take advantage of the discrete rotational
symmetry of the triple gaussian Hamiltonian.  Nevertheless, the tensor
decomposition provides sufficient compression of information to
facilitate efficient computation.

\subsubsection{Matrix Elements}
It is straightforward to calculate matrix elements for the Laplacian on
$R^1$:
\begin{equation}
  K_{mn} = \left\{\begin{array}{ll}
\frac{\hbar^2\pi^2}{\Delta{x}^2},&m=n\\
(-1)^{m-n}\cdot\frac{\hbar^2}{\Delta{x}^2\cdot(m-n)^2},& m\neq{n}\\
\end{array}\right.
\end{equation}
There is no closed form expression for the matrix elements of the
potential, but it is easy to perform numerical quadrature with these
functions.  For example, to compute
\begin{equation}
  V_{mn} = \int{G(x)\phi_m(x)\phi_n(x)\ dx}
\end{equation}
for $G(x)=e^{-\frac{x^2}{2\sigma^2}}$, one computes
\begin{equation}
  V_{mn} \approx
  \sum_{k=-N}^{N}{G(k\delta) \cdot \delta \cdot \phi_m(k\delta) \cdot
    \phi_n(k\delta)},
\end{equation}
where the stepsize $\delta$ should satisfy $\delta\leq\Delta{x}/2$.  It
is easy to show that the error is bounded by the sum of
\begin{equation}
  2\exp\of{-\frac{\abs{\sigma}^2\pi^2}{2\delta^2}},
\end{equation}
which controls the {\em aliasing error}, and
\begin{equation}
  \frac{\sqrt{2\pi\abs{\sigma}^2}}{\Delta{x}}
  \exp\of{-\frac{\of{N-1}^2\delta^2} {2\abs{\sigma}^2}}.
\end{equation}
which controls the {\em truncation error}.

\subsubsection{Other Program Parameters}
The grid spacing $\Delta{x}$ implies a limit on the maximum possible
momentum in a wave packet formed by this basis.  In order to obtain a
finite-rank operator, it is also necessary to limit the number of basis
functions.

The resonance computation used the following parameters:
\begin{enumerate}

  \item $X_0$, $X_1$, $Y_0$, and $Y_1$ are chosen to cover the region of
        the configuration space for which $V_3(x,y)\geq{10^{-4}}$.

  \item Let $L_x=X_1-X_0$ and $L_y=Y_1-Y_0$ denote the dimension of the
        computational domain.  The resonance calculation uses
        $N=N_x\cdot{N_y}$ basis functions, with
        $N_x=1.6\cdot\frac{L_x\sqrt{8E}}{2\pi\hbar}$ and
        $N_y=1.6\cdot\frac{L_y\sqrt{8E}}{2\pi\hbar}$.

  \item This gives
        \begin{equation}
           \begin{array}{lcl}
           \Delta{x}&=&L_x/N_x,\\
           \Delta{y}&=&L_y/N_y,\\
           \end{array}
        \end{equation}
        which limits the maximum momentum in a wave packet to
        $\abs{p_x}\leq\pi\hbar/\Delta{x}=1.6\sqrt{2E}$ and
        $\abs{p_y}\leq\pi\hbar/\Delta{y}=1.6\sqrt{2E}$.

\end{enumerate}

\section{Trapped Set Structure}
\subsection{\Poincare\ Section}
Because the phase space for the triple gaussian model is $R^4$ and its
flow is chaotic, a direct computation of the trapped set dimension is
difficult.  Instead, we try to compute its intersection with a {\em
{\Poincare\ }section}:

Let $E$ be a fixed energy, and recall that $R$ is the distance from each
gaussian bump to the origin.  Choose $R_0<R$ so that the circles $C_k$
of radius $R_0$ centered at each potential, for $k=0,1,2$, do not
intersect.  The angular momentum $p_\theta$ with respect to the $k$th
potential center is defined by
$p_\theta=\Delta{x}\cdot{p_y}-\Delta{y}\cdot{p_x}$, where
$\Delta{x}=x-R\cos\of{\theta_k}$, $\Delta{y}=y-R\sin\of{\theta_k}$, and
$\theta_k=\frac{2\pi k}{3}$.

\begin{figure}
  \begin{center}
  \PSbox{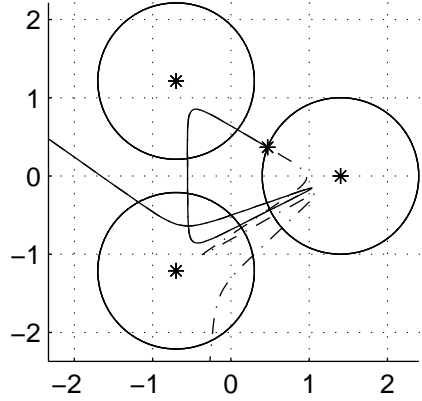}{2.5in}{2.5in}
  \end{center}

  \caption{A typical trajectory: Stars mark the potential centers.  In
  this case, $R=1.4$ and $E=0.5$.  The circles drawn in the figure have
  radius $1$, and the disjoint union of their cotangent bundles form the
  {\Poincare\ }section.  Trajectories start on the circle centered at
  bump 0 (the bumps are, counterclockwise, 0, 1, and 2) with some given
  angle $\theta$ and angular momentum $p_\theta$.  This trajectory
  generates the finite sequence $(\dot{0},1,2,0,2,\infty)$.  (Symbolic 
  sequences are discussed later in the paper.)  An illustration of
  resonance computation is depicted in Figure
  {\ref{fig:res-prog-params}}.  The dashed line is the time-reversed
  trajectory with the same initial conditions, generating the sequence
  $(\infty,2,0,2,\dot{0})$.}

  \label{fig:dim-prog-params}
\end{figure}

Let $P$ be the submanifold $P_0 \cup P_1 \cup P_2$ of $R^4$ (see Figure
{\ref{fig:dim-prog-params}}), where the coordinates $(\theta,p_\theta)$
in the submanifold $P_k$ are related to ambient phase space coordinates
$(x,y,p_x,p_y)$ by
\begin{equation}
  \label{eqn:poincare-embedding}
  \begin{array}{lcl}
  x&=&R\cos\of{\theta_k} + R_0\cos\of{\theta+\theta_k},\\
  y&=&R\sin\of{\theta_k} + R_0\sin\of{\theta+\theta_k},\\
  p_x&=&p_r\cos\of{\theta+\theta_k} -
        \frac{p_\theta}{R_0}\sin\of{\theta+\theta_k},\\
  p_y&=&p_r\sin\of{\theta+\theta_k} +
        \frac{p_\theta}{R_0}\cos\of{\theta+\theta_k}\\
  \end{array}
\end{equation}
and the radial momentum $p_r$ is
\begin{equation}
  p_r=\sqrt{E - V_3(x,y) - \frac{p_\theta^2}{2R_0^2}}.
\end{equation}
Note that this implicitly embeds $P$ into the energy surface
$\set{H=E}$, and the radial momentum $p_r$ is always positive: The
vector $(p_x,p_y)$ points away from the center of $C_k$.

The trapped set is naturally partitioned into two subsets: The first
consists of trajectories which visit all three bumps, the second of
trajectories which bounce between two bumps.  The second set forms a
one-dimensional subspace of $K_E$, so the {\em finite stability} of the
Minkowski dimension\footnote{That is, $D(A\cup{B})=\max\of{D(A),D(B)}$.
For details, see {\cite{falconer}}.} implies that the second set does
not contribute to the dimension of the trapped set.  More importantly,
most trajectories which visit all three bumps will also cut through $P$.

One can thus reduce the dimension of the problem by restricting the flow
to $K_E\cap{P}$, as follows: Take any point $(\theta,p_\theta)$ in
$P_k$, and form the corresponding point $(x,y,p_x,p_y)$ in $R^4$ via
Equation (\ref{eqn:poincare-embedding}).  Follow along the trajectory
$\Phi_t(x,y,p_x,p_y)$.  If the trajectory does not escape, eventually it
must encounter one of the other circles, say $C_{k'}$.  Generically,
trajectories cross $C_{k'}$ twice at each encounter, and we denote the
coordinates $(\theta',p_\theta')$ (in $P_{k'}$) of the {\em outgoing
intersection} by
\begin{equation}
  \label{eqn:PhiP-def}
  \PhiP(\theta,p_\theta,k)=(\theta',p_\theta',k').
\end{equation}
If a trajectory escapes from the trapping region, we can symbolically
assign $\infty$ to $\PhiP$.  The map $\PhiP$ then generates stroboscopic
recordings of the flow $\Phi_t$ on the submanifold $P$, and the
corresponding discrete dynamical system has trapped set $K_E\cap{P}$.
So, instead of computing $\Phi_t$ on $R^4$, one only needs to compute
$\PhiP$ on $P$.  By symmetry, it will suffice to compute the dimension
of $\KEhat=K_E\cap{P_0}$.  Pushing $\KEhat$ forward along the flow
$\Phi_t$ adds one dimension, so $D\of{K_E}=D\of{\KEhat}+1$.  Being a
subset of the two-dimensional space $P_0$, $\KEhat$ is easier to work
with.

Readers interested in a more detailed discussion of {\Poincare\
}sections and their use in dynamics are referred to {\cite{gjs}}.  For
an application to the similar but simpler setting of hard disc
scattering, see {\cite{GR1,gaspard}}.  Also, Knauf has applied some of
these ideas in a theoretical investigation of classical scattering by
Coulombic potentials {\cite{knauf}}.

\subsection{Self-Similarity}
\begin{figure}
  \begin{center}
  \PSbox{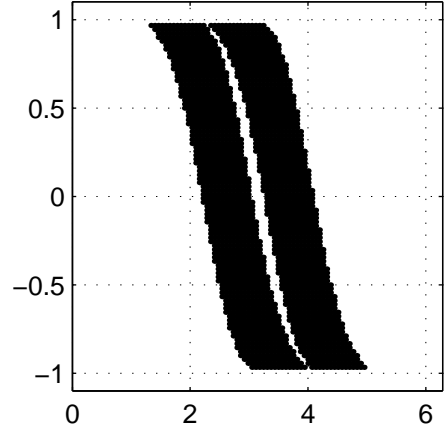}{2.5in}{2.5in}
  \end{center}

  \caption{Points in $P_0$ which do not go to $\infty$ after one
  iteration of $\PhiP$.  The horizontal axis is $\theta$ and the
  vertical axis is $p_\theta$.}

  \label{fig:bounce1}
\end{figure}

Much is known about the self-similar structure of the trapped set for
hard disc scattering {\cite{GR1,gaspard}}; less is known about ``soft
scatterers'' like the triple gaussian system.  However, computational
results and analogy with hard disc scattering give strong support to the
idea that $K$ (and hence $\KEhat$) is self-similar.\footnote{More
precisely, {\em self-affine}.} Consider Figures {\ref{fig:bounce1}} -
{\ref{fig:bounce7}}: They show clearly that $\KEhat$ is self-similar.
(In these images, $E=0.5$ and $R_0=1.0$.)  However, it is also clear
that, unlike objects such as the Cantor set or the Sierpi\'{n}ski
gasket, $\KEhat$ is not exactly self-similar.
\begin{figure}
  \begin{center}
  \PSbox{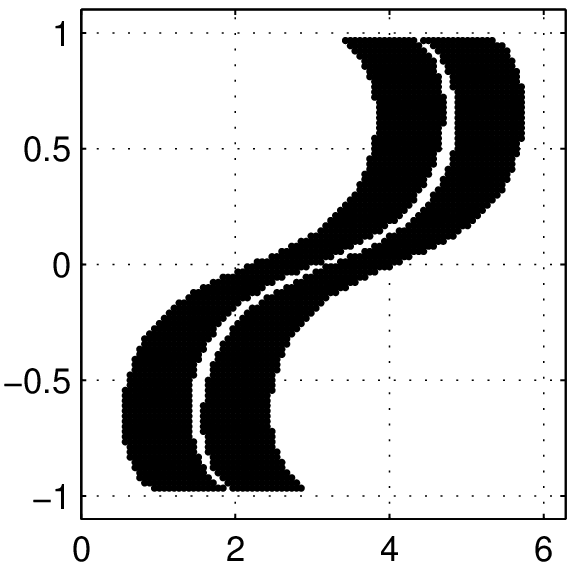}{2.5in}{2.5in}
  \end{center}

  \caption{Points in $P_0$ which do not go to $\infty$ after one
  iteration of $\PhiP^{-1}$.  The horizontal axis is $\theta$ and the
  vertical axis is $p_\theta$.}

  \label{fig:bounce2}
\end{figure}
\begin{figure}
  \begin{center}
  \PSbox{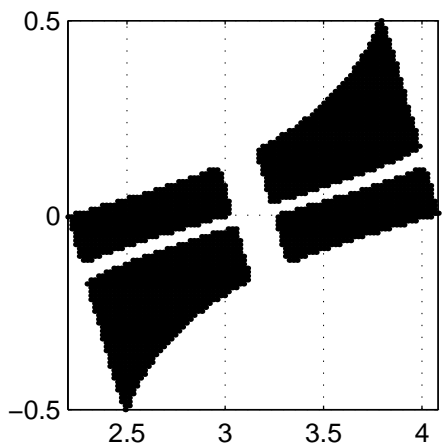}{2.5in}{2.5in}
  \end{center}

  \caption{The intersection of the sets in Figures {\ref{fig:bounce1}}
  and \ref{fig:bounce2}.  These points correspond to symmetric sequences
  of length 3.}

  \label{fig:bounce3}
\end{figure}
\begin{figure}
  \begin{center}
  \PSbox{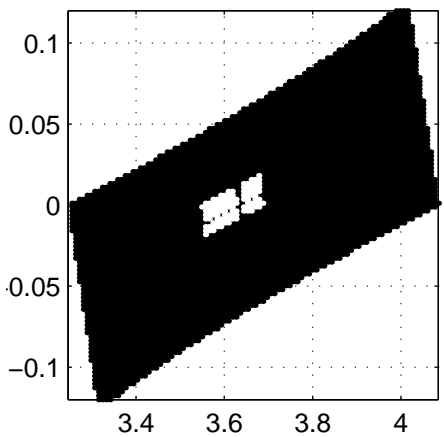}{2.5in}{2.5in}
  \end{center}

  \caption{The lower-right ``island'' in Figure {\ref{fig:bounce3}},
  magnified.  The white cut-out in the middle is the subset
  corresponding to symmetric sequences of length 5.}

  \label{fig:bounce4}
\end{figure}
\begin{figure}
  \begin{center}
  \PSbox{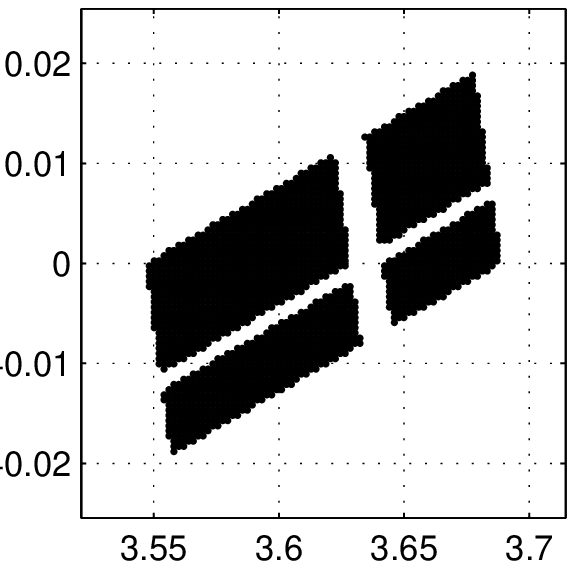}{2.5in}{2.5in}
  \end{center}

  \caption{The cut-out part of Figure {\ref{fig:bounce4}}, magnified.
  Recall that these correspond to symmetric sequences of length 5;
  compare with Figure {\ref{fig:bounce3}}.}

  \label{fig:bounce5}
\end{figure}
\begin{figure}
  \begin{center}
  \PSbox{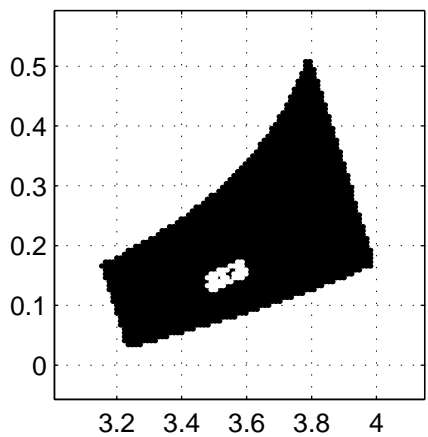}{2.5in}{2.5in}
  \end{center}

  \caption{The upper-right island in Figure {\ref{fig:bounce3}}.  The
  white cut-out in the middle is, again, the subset corresponding to
  symmetric sequences of length 5.}

  \label{fig:bounce6}
\end{figure}
\begin{figure}
  \begin{center}
  \PSbox{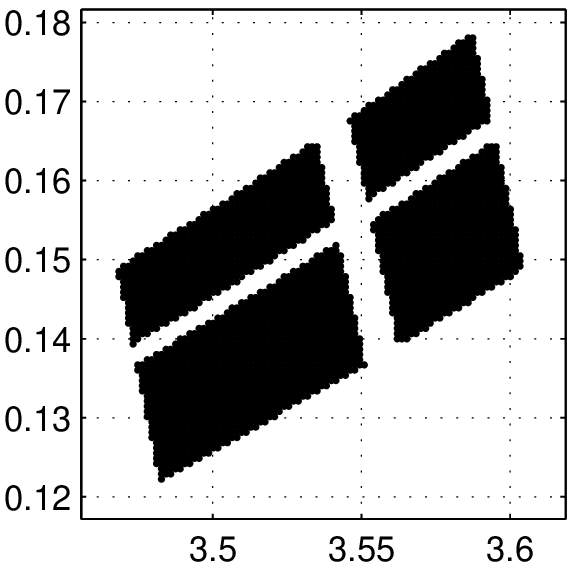}{2.5in}{2.5in}
  \end{center}

  \caption{The cut-out part of Figure {\ref{fig:bounce6}}, magnified.
  Recall that these correspond to symmetric sequences of length 5.
  Compare with Figures {\ref{fig:bounce3}} and {\ref{fig:bounce5}}.}

  \label{fig:bounce7}
\end{figure}

\subsection{Symbolic Dynamics}
The computation of $D\of{\KEhat}$ uses symbolic sequences, which
requires a brief explanation: For any point $(\theta,p_\theta)$, let
$s_i$ denote the third component of $\PhiP^i(\theta,p_\theta)$ (see
(\ref{eqn:PhiP-def})), for any integer $i$.  Thus, $s_i$ is the index
$k$ of the circle $C_k$ that the trajectory intersects at the $i$th
iteration of $\PhiP$ (or the $\abs{i}$th iteration of $\PhiP^{-1}$).
Such symbolic sequences $s=(...s_i,s_{i+1}...)$ satisfy
$s_i\in\set{0,1,2,\infty}$ and $s_i\neq{s_{i+1}}$ for all $i$, and with
$\infty$ occuring only at the ends.  Let us call sequences satisfying
these conditions {\em valid}.

For example, the trajectory in Figure {\ref{fig:res-prog-params}}
generates the valid sequence $(\dot{0},1,2,0,2,\infty)$, where the dot
over $0$ indicates that the initial point $(\theta,p_\theta)$ of the
trajectory belongs to $P_0$.  Thus, we can label collections of
trajectories using valid sequences, and label points in $P$ with
``dotted'' sequences.  Clearly, trapped trajectories generate
bi-infinite sequences.\footnote{In hard disc scattering, the converse
holds for sufficiently large $R$: To each bi-infinite valid sequence
there exists a trapped trajectory generating that sequence.  This may
not hold in the triple gaussian model, and in any case it is not
necessary for the computation.}

The islands in Figures {\ref{fig:bounce3}} - {\ref{fig:bounce6}}
correspond to {\em symmetric sequences centered at 0}, of the form
$s=(...s_{-k},...,s_{-1},\dot{0},s_{1},s_{2},...,s_{k}...)$: By keeping
track of the symbolic sequences generated by each trajectory, one can
easily label and isolate each island.  This is a useful property from
the computational point of view.

\subsection{Dimension Estimates}
To compute the Minkowski dimension using Equation
(\ref{eqn:minkowski-def1}), we need to determine when a given point is
within $\epsilon$ of $\KEhat$.  This is generally impossible: The best
one can do is to generate longer and longer trajectories which stay
trapped for increasing (but finite) amounts of time.

Instead, one can estimate a closely related quantity, the {\em
information dimension}, in the following way: Let $\KEhatt{k}$ denote
the set of all points in $P_0$ corresponding to symmetric sequences of
length $2k+1$ centered at 0.  That is, $\KEhatt{k}$ consists of all
points in $P_0$ which generate trajectories (both forwards and backwards
in time) that bounce at least $k$ times before escaping.  The sets
$\KEhatt{k}$ decrease monotonically to $\KEhat$:
$\KEhatt{k}\supset\KEhatt{k+1}$ and
$\cap_{k=0}^\infty{\KEhatt{k}}=\KEhat$.

One can then estimate the information dimension using the following
algorithm:
\begin{enumerate}

  \item {\em Initialization:} Cover $P_0$ with a mesh $L_0$ with
  $N_0\times{N_0}$ grid points and mesh size $\epsilon_0$.

  \item {\em Recursion:} Begin with $\KEhatt{1}$, which consists of four
  islands corresponding to symmetric sequences of length $2\cdot{1}+1=3$
  (see Figure {\ref{fig:bounce3}}).  Magnify each of these islands and
  compute the sub-islands corresponding to symmetric sequences of length
  $5$ (see Figures {\ref{fig:bounce4}} and {\ref{fig:bounce6}}).  Repeat
  this procedure to recursively compute the islands of $\KEhatt{k+1}$
  from those of $\KEhatt{k}$.  Continue until $k=k_0$, where $k_0$ is
  sufficiently large that each island of $\KEhatt{k_0}$ has diameter
  smaller than the mesh size $\epsilon_0$ of $L_0$.

  \item {\em Estimation:} Using the islands of $\KEhatt{k_0}$, estimate
  the probability
  \begin{equation}
    p_{ij}=\frac{\vol\of{\KEhatt{k_0} \cap B_{ij}}}
           {\vol\of{\KEhatt{k_0}}}
  \end{equation}
  for the $(ij)$th cell of $L_0$.  We can then compute the dimension via
  \begin{equation}
    \label{eqn:infodim}
    D\of{\KEhat} \approx
                 \frac{-\sum_{ij}{p_{ij}\log\of{p_{ij}}}}{\log\of{N_0}},
  \end{equation}
  which reduces to (\ref{eqn:minkowski-def1}) when the distribution is
  uniform because $\epsilon_0\sim{1/N_0}$.

\end{enumerate}
\begin{figure}
  \begin{center}
  \PSbox{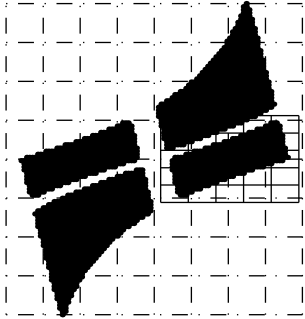}{2.5in}{2.5in}
  \end{center}

  \caption{This figure illustrates the recursive step in the dimension
  estimation algorithm: The dashed lines represent $L_0$, while the
  solid lines represent a smaller mesh centered on one of the islands.
  The $N_0\times{N_0}$ mesh $L_0$ remains fixed throughout the
  computation, but the smaller $N_1\times{N_1}$ mesh is constructed for
  each island of $\KEhatt{k}$ up to the value of $k$ specified by the
  algorithm.}

  \label{fig:dimzoom}
\end{figure}
Under suitable conditions (as is assumed to be the case here), the
information dimension agrees with both the Hausdorff and the Minkowski
dimensions.\footnote{See {\cite{pesin}} for a discussion of the
relationship between these dimensions, as well as their use in
multifractal theory.}

The algorithm begins with the lattice $L_0$ with which one wishes to
compute the dimension.  It then recursively computes $\KEhatt{k}$ for
for increasing values of $k$, until it closely approximates $\KEhat$
relative to the mesh size of $L_0$.  It is easy to keep track of points
belonging to each island in this computation, since each island
corresponds uniquely to a finite symmetric sequence.  Note that while
the large mesh $L_0$ remains fixed throughout the computation, the
recursive steps require smaller $N_1\times{N_1}$ meshes around each
island of $\KEhatt{k}$ up to the value of $k$ specified by the
algorithm.  See Figure {\ref{fig:dimzoom}}.

\section{Numerical Results}
\subsection{Resonance Counting}
\begin{figure}
  \begin{center}
  \PSbox{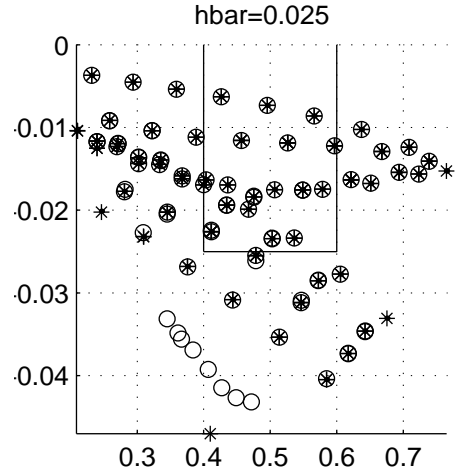}{2.5in}{2.5in}
  \end{center}

  \caption{These are the eigenvalues of $\Hnalpha$, for $E=0.5$,
  $\hbar=0.025$, $R=1.4$, and $\alpha\in\set{0.0624, 0.0799, 0.0973}$.
  This calculation used an $102\times{108}$ grid, and $90$ out of
  $N=11016$ eigenvalues were computed.}

  \label{fig:resplot0}
\end{figure}
\begin{figure}
  \begin{center}
  \PSbox{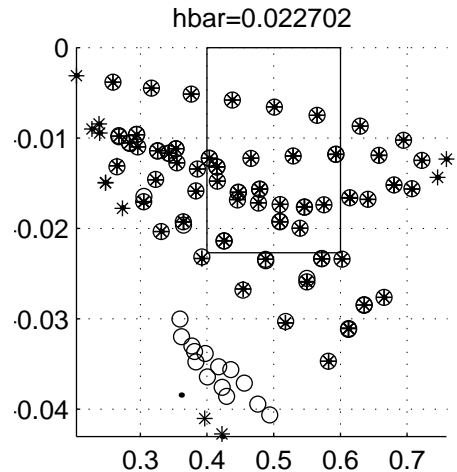}{2.5in}{2.5in}
  \end{center}

  \caption{Eigenvalues for $E=0.5$, $\hbar=0.022702$, $R=1.4$, and
  $\alpha\in\set{0.0567, 0.0741, 0.0916}$, using $112\times{119}$ grid
  and $98$ out of $N=13328$ eigenvalues.}

  \label{fig:resplot1}
\end{figure}
\begin{figure}
  \begin{center}
  \PSbox{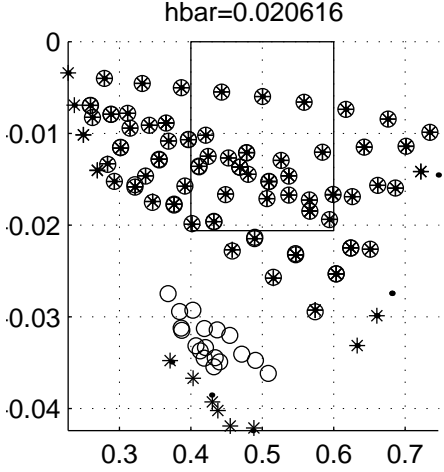}{2.5in}{2.5in}
  \end{center}

  \caption{Eigenvalues for $E=0.5$, $\hbar=0.020616$, $R=1.4$, and
  $\alpha\in\set{0.0515, 0.0689, 0.0864}$, using $123\times{131}$ grid
  and $107$ out of $N=16113$ eigenvalues.}

  \label{fig:resplot2}
\end{figure}
\begin{figure}
  \begin{center}
  \PSbox{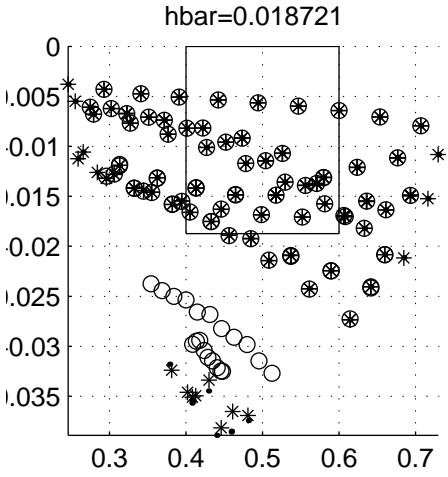}{2.5in}{2.5in}
  \end{center}

  \caption{Eigenvalues for $E=0.5$, $\hbar=0.018721$, $R=1.4$, and
  $\alpha\in\set{0.0468, 0.0642, 0.0817}$, using $135\times{144}$ grid
  and $116$ out of $N=19440$.}

  \label{fig:resplot3}
\end{figure}
\begin{figure}
  \begin{center}
  \PSbox{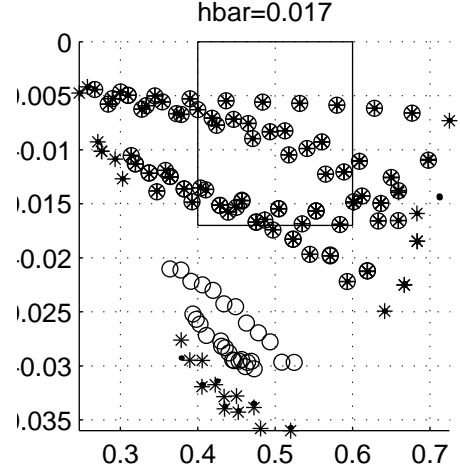}{2.5in}{2.5in}
  \end{center}

  \caption{Eigenvalues for $E=0.5$, $\hbar=0.017$, $R=1.4$, and
  $\alpha\in\set{0.0425, 0.0599, 0.0774}$, using $149\times{159}$ grid
  and $127$ out of $N=23691$ eigenvalues.}

  \label{fig:resplot4}
\end{figure}
As an illustration of complex scaling, Figures {\ref{fig:resplot0}} -
{\ref{fig:resplot4}} contain resonances for $R=1.4$ and
$\hbar\in\interval{0.017,0.025}$.  Eigenvalues of $\Hnalpha$ for
different values of $\alpha$ are marked by different styles of points,
and the box has depth $\hbar$ and width $0.2$, with $E_0=0.4$ and
$E_1=0.6$.  These plots may seem somewhat empty because only those
eigenvalues of $\Hnalpha$ in regions of interest were computed.  Notice
the cluster of eigenvalues near the bottom edge of the plots: These are
{\em not} resonances because they vary under perturbations in $\alpha$.
Instead, they belong to an approximation of the (scaled) continuous
spectrum.

\begin{figure}
  \begin{center}
  \PSbox{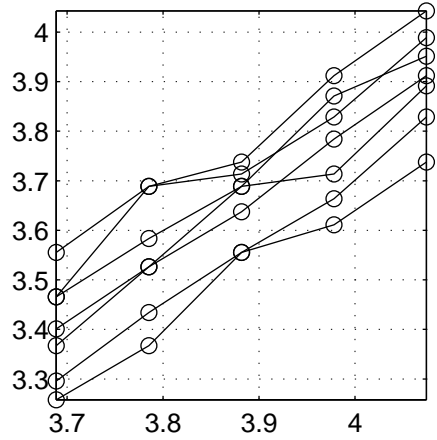}{2.5in}{2.5in}
  \end{center}

  \caption{$\log\of{N_{res}}$ as a function of $-\log\of{\hbar}$, for
  $\hbar$ varying from $0.017$ to $0.025$ and $R=1.4+0.05\cdot{k}$, with
  $0\leq{k}\leq{6}$.  (The lowest curve corresponds to $R=1.7$, while
  the highest curve corresponds to $R=1.4$.)}

  \label{fig:resfun2}
\end{figure}
\begin{figure}
  \begin{center}
  \PSbox{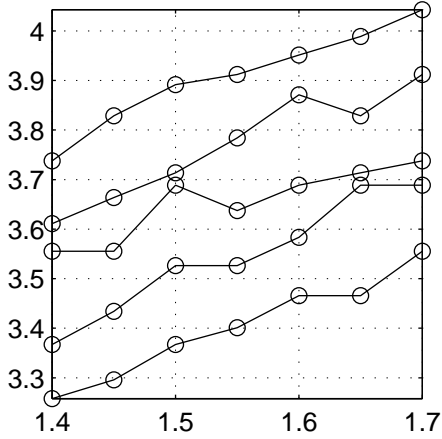}{2.5in}{2.5in}
  \end{center}

  \caption{$\log\of{N_{res}}$ as a function of $R$, for different values
  of $-\log\of{\hbar}$: The highest curve corresponds to $\hbar=0.017$,
  while the lowest curve corresponds to $\hbar=0.025$.}

  \label{fig:resfun3}
\end{figure}
\begin{figure}
  \begin{center}
  \PSbox{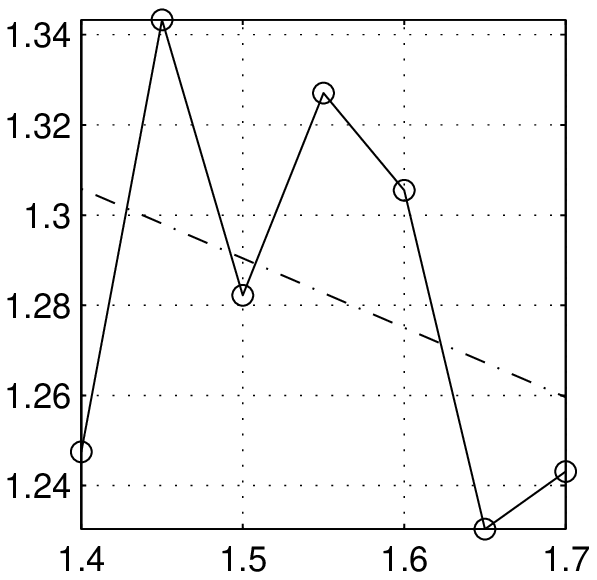}{2.5in}{2.5in}
  \end{center}

  \caption{The slopes extracted from Figure {\ref{fig:resfun2}}, as a
  function of $R$.  The dotted curve is a least-squares regression of
  the ``noisy'' curve.}

  \label{fig:resfun4}
\end{figure}
It is more interesting to see $\log\of{N_{res}}$ as a function of
$-\log\of{\hbar}$ and $R$.  This is shown in Figures {\ref{fig:resfun2}}
and {\ref{fig:resfun3}}.  Using least-squares regression, we can extract
approximate slopes for the curves in Figure {\ref{fig:resfun2}}; these
are shown in Table {\ref{tab:compare}} and plotted in Figure
{\ref{fig:resfun4}}.

\subsection{Trapped Set Dimension}
\begin{figure}
  \begin{center}
  \PSbox{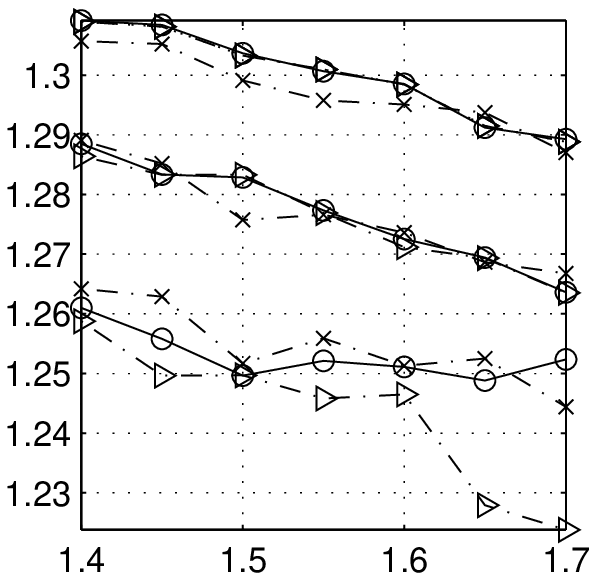}{2.5in}{2.5in}
  \end{center}

  \caption{This figure shows $\frac{D\of{K_E}+1}{2}$ as a function of
  $R$: The top group of curves have $E=0.4$, the middle $E=0.5$, and the
  bottom $E=0.6$.  Solid curves marked with circles represent
  computations where $N_0=10000$, $N_1=100$,
  $\frac{2\pi}{3}\leq\theta\leq\frac{4\pi}{3}$, and
  $-\frac{1}{2}\leq{p_\theta}\leq\frac{1}{2}$.  Dashed curves marked
  with X's represent computations where $N_0=14142$, whereas dashed
  curves marked with triangles represent computations where $N_0=10000$
  and $N_1=71$.  The recursion depth $k_0$ in all these figues is $6$.
  The $E=0.6$ curve does not appear to have completely converged but
  suffices for our purpose here.}

  \label{fig:dimconv}
\end{figure}
\begin{table}
  \begin{center}
  \begin{tabular}{|c|c|c|c|}
    \hline
     $R$  & $E=0.4$  & $E=0.5$&$E=0.6$\\\hline\hline
     1.4 & 1.3092 & 1.2885 & 1.261\\\hline
     1.45 & 1.3084 & 1.2834 & 1.2558\\\hline
     1.5 & 1.3037 & 1.2829 & 1.2497\\\hline
     1.55 & 1.3007 & 1.2773 & 1.2521\\\hline
     1.6 & 1.2986 & 1.2725 & 1.2511\\\hline
     1.65 & 1.2912 & 1.2694 & 1.2488\\\hline
     1.7 & 1.2893 & 1.2636 & 1.2524\\\hline
  \end{tabular}
  \end{center}

  \caption{Estimates of $\frac{D\of{K_E}+1}{2}$ as a function of $R$.}

  \label{tab:dims}
\end{table}
For comparison, $\frac{D\of{K_E}+1}{2}$ is plotted as a function of $R$
in Figure {\ref{fig:dimconv}}.  The figure contains curves corresponding
to different energies $E$: The top curve corresponds to $E=0.4$, the
middle curve $E=0.5$, and the bottom curve $E=0.6$.  It also contains
curves corresponding to different program parameters, to test the
numerical convergence of dimension estimates.  These curves were
computed with $\theta\in\interval{\frac{2\pi}{3},\frac{4\pi}{3}}$,
$p_\theta\in\interval{-\frac{1}{2},\frac{1}{2}}$, and recursion depth
$k_0=6$ (corresponding to symmetric sequences of length $2\cdot{6}+1
=13$); the caption contains the values of $N_0$ and $N_1$ for each
curve.  For reference, Table {\ref{tab:dims}} contains the dimension
estimates shown in the graph.  It is important to note that, while the
dimension does depend on $E$ and $R$, it only does so weakly: Relative
to its value, $\frac{D\of{K_E}+1}{2}$ is very roughly constant across
the range of $R$ and $E$ computed here.

\subsection{Discussion}
\begin{table}
  \begin{center}
  \begin{tabular}{|c|c|c|c|}
    \hline
    $R$ & slope   & $\frac{D\of{K_E}+1}{2}$ & relative error \\\hline\hline
   1.4  & 1.2475  & 1.2885   & 0.032888               \\\hline
   1.45 & 1.3433  & 1.2834   & 0.044645               \\\hline
   1.5  & 1.2822  & 1.2829   & 0.00052244             \\\hline
   1.55 & 1.327   & 1.2773   & 0.037472               \\\hline
   1.6  & 1.3055  & 1.2725   & 0.025256               \\\hline
   1.65 & 1.2304  & 1.2694   & 0.031756               \\\hline
   1.7  & 1.2431  & 1.2636   & 0.016509               \\\hline
  \end{tabular}
  \end{center}

  \caption{This table shows the slopes extracted from Figure
  {\ref{fig:resfun2}}, as well as the scaling exponents one would expect
  if the conjecture were true ($\of{D\of{K_E}+1}/2$, computed at
  $E=0.5$).  Relative errors are also shown.}

  \label{tab:compare}
\end{table}
\begin{figure}
  \begin{center}
  \PSbox{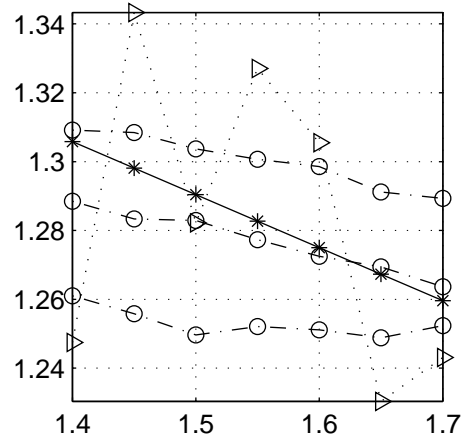}{2.5in}{2.5in}
  \end{center}

  \caption{Dashed lines with circles represent $\frac{D\of{K_E}+1}{2}$
  as functions of $R$, for $E\in\set{0.4,0.5,0.6}$.  The dotted curve
  with triangles is the scaling exponent curve from Figure
  {\ref{fig:resfun4}}, while the solid curve with stars is the linear
  regression curve from that figure.  Relative to the value of the
  dimension, the fluctuations are actually fairly small: See Table
  {\ref{tab:compare}} for a quantitative comparison.}

  \label{fig:compare}
\end{figure}
Table {\ref{tab:compare}} contains a comparison of
$\frac{D\of{K_E}+1}{2}$ (for $E=0.5$) as a function of $R$, versus the
scaling exponents from Figure {\ref{fig:resfun4}}.  Figure
{\ref{fig:compare}} is a graphical representation of similar
information.  This figure shows that even though the scaling curve in
Figure {\ref{fig:resfun4}} is noisy, its trend nevertheless agrees with
the conjecture.  Furthermore, the relative size of the fluctuations is
small.  At the present time, the source of the fluctuation is not known,
but it is possibly due to the fact that the range of $\hbar$ explored
here is simply too large to exhibit the asymptotic behavior
clearly.\footnote{But see Footnote {\ref{fn:quad}}.}

\begin{figure}
  \begin{center}
  \PSbox{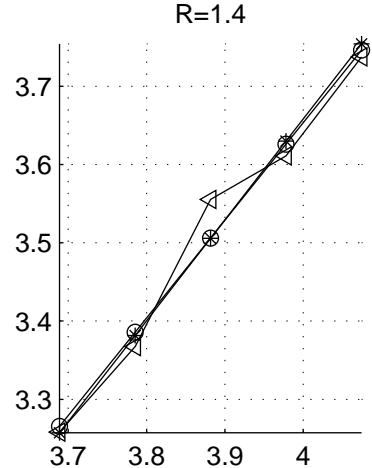}{2.5in}{2.5in}
  \end{center}

  \caption{For $R=1.4$: Triangles represent numerical data, circles
  least squares regression, and stars the slope predicted by the
  conjecture.  $\hbar$ ranges from $0.025$ down to $0.017$.}

  \label{fig:hist0}
\end{figure}

\begin{figure}
  \begin{center}
  \PSbox{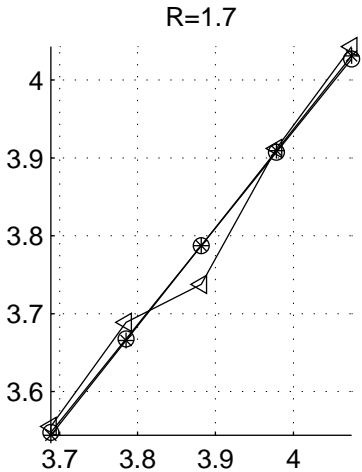}{2.5in}{2.5in}
  \end{center}
  \caption{$R=1.7$, $0.017\leq{\hbar}\leq{0.025}$.}
  \label{fig:hist6}
\end{figure}

Figures {\ref{fig:hist0}} - {\ref{fig:hist6}} contain plots of
$\log\of{N_{res}}$ versus $-\log\of{\hbar}$, for various values of $R$.
Along with the numerical data, the least-squares linear fit and the
scaling law predicted by the conjecture are also plotted.\footnote{The
conjecture only supplies the exponents for power laws, not the constant
factors.  In the context of these logarithmic plots, this means the
conjecture gives us only the slopes, not the vertical shifts.  It was
thus necessary to compute an $y$-intercept for each ``prediction'' curve
(for the scaling law predicted by the conjecture) using least squares.}
In contrast with Figure {\ref{fig:compare}}, these show clear agrement
between the asymptotic distribution of resonances and the scaling
exponent predicted by the conjecture.
\begin{figure}
  \begin{center}
  \PSbox{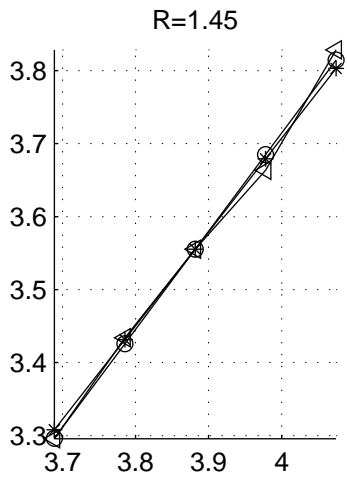}{2.5in}{2.5in}
  \end{center}

  \caption{Same for $R=1.45$.  Again, $\hbar$ ranges from $0.025$ to
  $0.017$.}

  \label{fig:hist1}
\end{figure}
\begin{figure}
  \begin{center}
  \PSbox{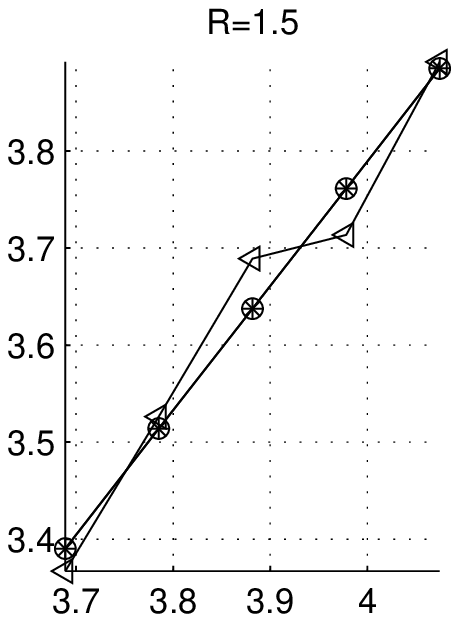}{2.5in}{2.5in}
  \end{center}
  \caption{$R=1.5$, $0.017\leq{\hbar}\leq{0.025}$.}
  \label{fig:hist2}
\end{figure}
\begin{figure}
  \begin{center}
  \PSbox{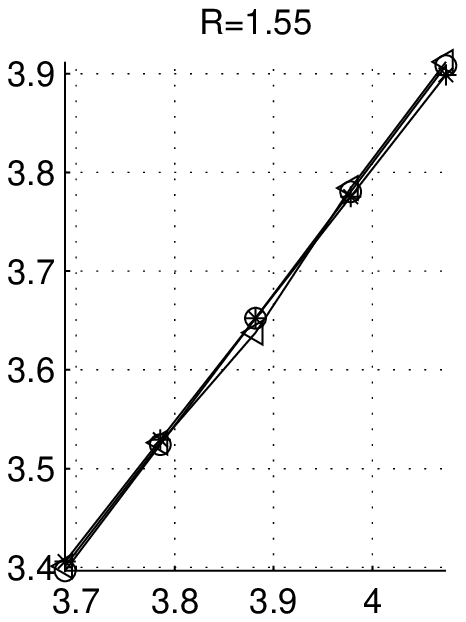}{2.5in}{2.5in}
  \end{center}
  \caption{$R=1.55$}
  \label{fig:hist3}
\end{figure}
\begin{figure}
  \begin{center}
  \PSbox{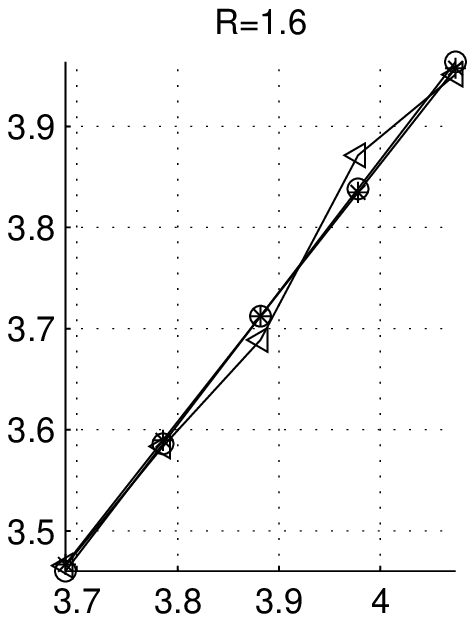}{2.5in}{2.5in}
  \end{center}
  \caption{$R=1.6$, $0.017\leq{\hbar}\leq{0.025}$.}
  \label{fig:hist4}
\end{figure}
\begin{figure}
  \begin{center}
  \PSbox{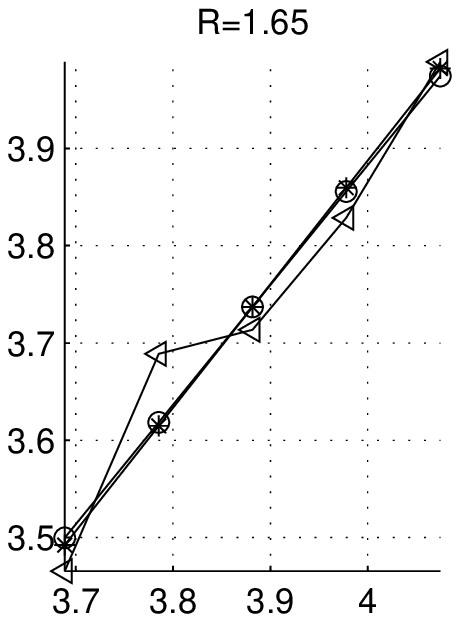}{2.5in}{2.5in}
  \end{center}
  \caption{$R=1.65$, $0.017\leq{\hbar}\leq{0.025}$.}
  \label{fig:hist5}
\end{figure}

\subsection{Double Gaussian Scattering}
Finally, we compute resonances for the double gaussian model (setting
$m=2$ in (\ref{eqn:potdef}).  This case is interesting for two reasons:
First, there exist rigorous results {\cite{GS,miller0}} against which we
can check the correctness of our results.  Second, it helps determine
the validity of semiclassical arguments for the values of $\hbar$ used
in computing resonances for the triple gaussian model.

\begin{figure}
  \begin{center}
  \PSbox{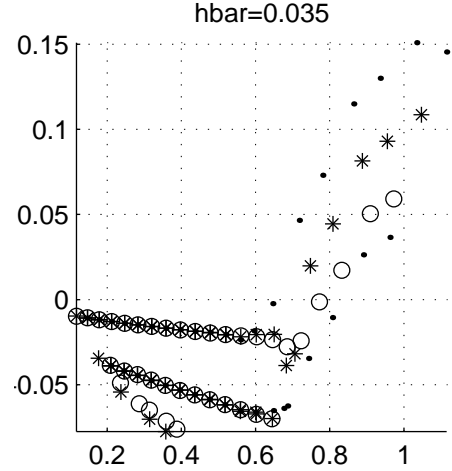}{2.5in}{2.5in}
  \end{center}
  \caption{Resonances for two-bump scattering with $\hbar=0.035$.}
  \label{fig:twobumps0}
\end{figure}

\begin{figure}
  \begin{center}
  \PSbox{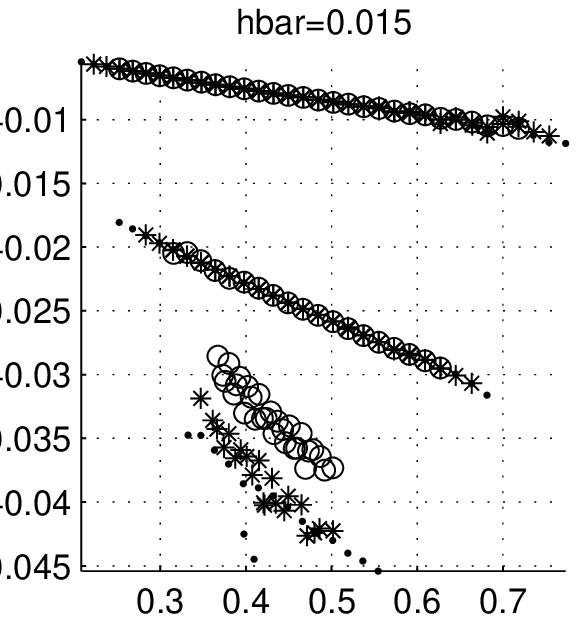}{2.5in}{2.5in}
  \end{center}
  \caption{Resonances for two-bump scattering with $\hbar=0.015$.}
  \label{fig:twobumps6}
\end{figure}

The resonances are shown in Figures {\ref{fig:twobumps0}} -
{\ref{fig:twobumps6}}: In these plots, $R=1.4$ and $\hbar$ ranges from
$0.035$ to $0.015$.  One can observe apparent pseudospectral effects in
the first few figures {\cite{trefethen,z3}}; this is most likely because
the scaling angle $\alpha$ used here is twice as large as suggested in
Section {\ref{sec:alpha-choice}}, to exhibit the structure of resonances
farther away from the real axis.
\begin{figure}
  \begin{center}
  \PSbox{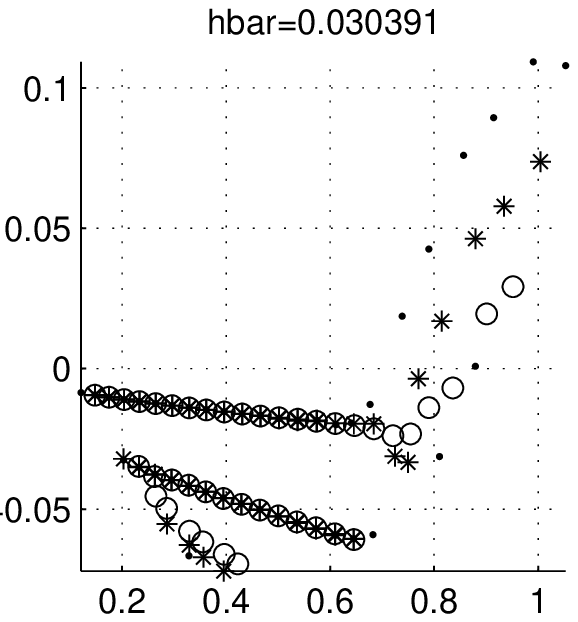}{2.5in}{2.5in}
  \end{center}
  \caption{Resonances for two-bump scattering with $\hbar=0.030391$.}
  \label{fig:twobumps1}
\end{figure}
\begin{figure}
  \begin{center}
  \PSbox{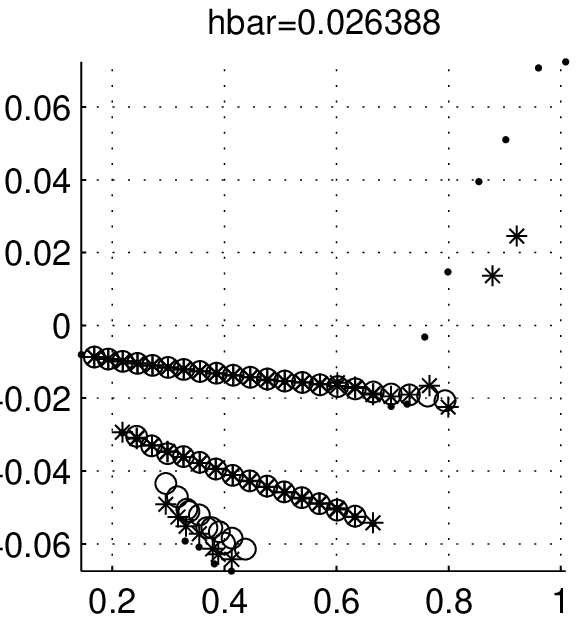}{2.5in}{2.5in}
  \end{center}
  \caption{Resonances for two-bump scattering with $\hbar=0.026388$.}
  \label{fig:twobumps2}
\end{figure}
\begin{figure}
  \begin{center}
  \PSbox{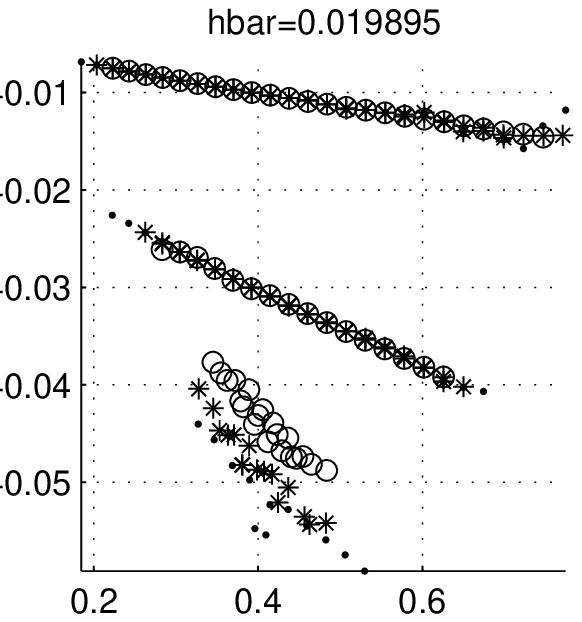}{2.5in}{2.5in}
  \end{center}
  \caption{Resonances for two-bump scattering with $\hbar=0.019895$.}
  \label{fig:twobumps4}
\end{figure}
\begin{figure}
  \begin{center}
  \PSbox{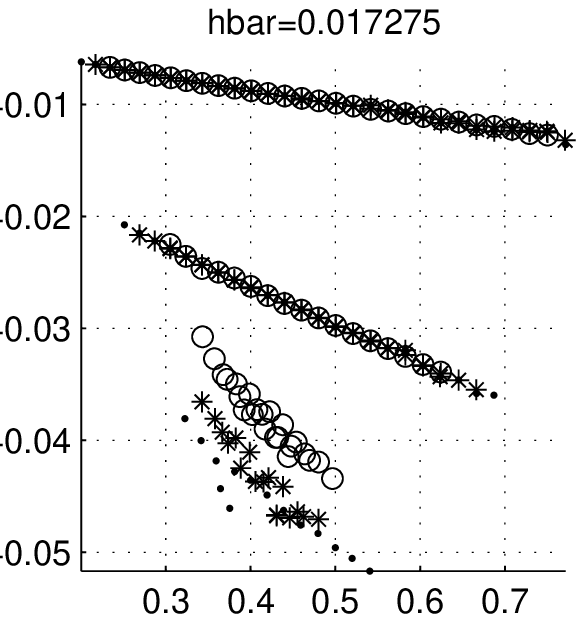}{2.5in}{2.5in}
  \end{center}
  \caption{Resonances for two-bump scattering with $\hbar=0.017275$.}
  \label{fig:twobumps5}
\end{figure}

To compare this information with known results {\cite{GS,miller0}}, we
need some definitions: For a given energy $0<E<1$, define $C(E)$ by
\begin{equation}
  C(E)=2\int_{x_0(E)}^{x_1(E)} {\sqrt{2\cdot\of{E - V\of{x}}} dx},
\end{equation}
where the limits of integration are
\begin{equation}
  \begin{array}{lcl}
    x_0(E)&=&-R+\sqrt{-2\sigma^2\log\of{E}},\\
    x_1(E)&=& R-\sqrt{-2\sigma^2\log\of{E}}.\\
  \end{array}
\end{equation}
Let $\theta(E)$ denote the larger (in absolute value) eigenvalue of
$D\PhiP(0,0)$; $\log\of{\theta}$ is the Lyapunov exponent of $\PhiP$,
and is easy to compute numerically in this case.  Note that for two-bump
scattering, each energy $E$ determines a unique periodic trapped
trajectory, and $C(E)$ is the classical action computed along that
trajectory.

\begin{figure}
  \begin{center}
  \PSbox{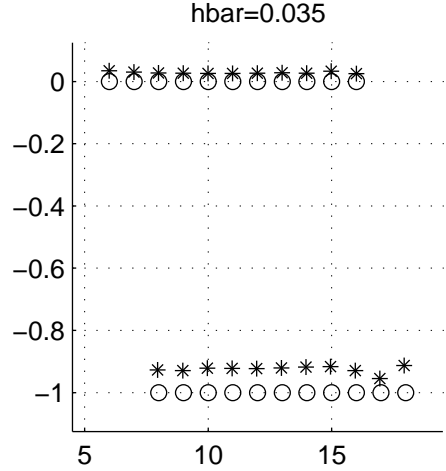}{2.5in}{2.5in}
  \end{center}
  \caption{Lattice points for $\hbar=0.035$.}
  \label{fig:latt0}
\end{figure}

\begin{figure}
  \begin{center}
  \PSbox{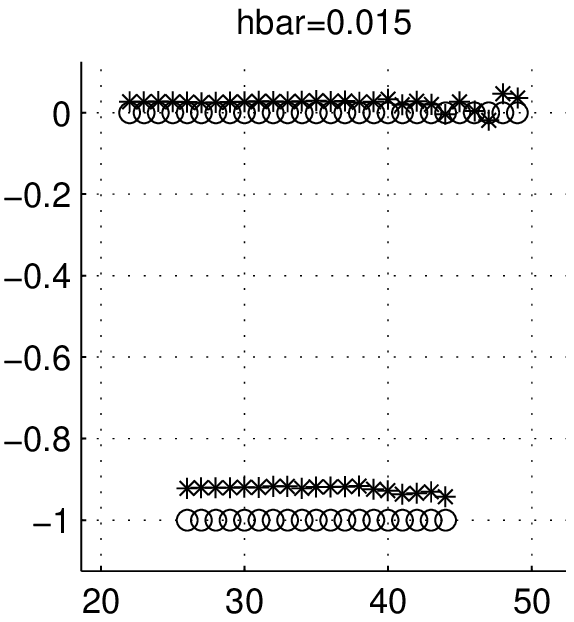}{2.5in}{2.5in}
  \end{center}
  \caption{Lattice points for $\hbar=0.015$.}
  \label{fig:latt6}
\end{figure}

Since these expressions are analytic, they have continuations to a
neighborhood of the real line --- $C(E)$ becomes a contour integral.  In
{\cite{GS}}, it was shown that any resonance $\lambda=E-i\gamma$ must
satisfy
\begin{equation}
  \label{eqn:GS}
  \begin{array}{lcl}
    C(\lambda)&=&2\pi\hbar\of{m+\frac{1}{2}} -\\
      &&i\hbar\of{n+\frac{1}{2}}\log\of{\theta\of{\real{\lambda}}} +\\
      &&O\of{\hbar^2},\\
  \end{array}
\end{equation}
where $m$ and $n$ are nonnegative integers.  (The $\frac{1}{2}$ in
$m+\frac{1}{2}$ comes from the Maslov index associated with the
classical turning points.)  This suggests that we define the map
$F\of{\lambda}=\of{F_1\of{\lambda},F_2\of{\lambda}}$, where
\begin{equation}
  F_1\of{\lambda} = \frac{\real{C(\lambda)}}{2\pi\hbar}-\frac{1}{2}
\end{equation}
and
\begin{equation}
  F_2\of{\lambda} = \frac{\imag{C(\lambda)}}
                    {\hbar\log\of{\theta\of{\real{\lambda}}}} +
                    \frac{1}{2}.
\end{equation}
$F$ should map resonances to points on the square integer lattice, and
this is indeed the case: Figures {\ref{fig:latt0}} - {\ref{fig:latt6}}
contain images of resonances under $F$, with circles marking the nearest
lattice points.  The agreement is quite good, in view of the fact that
we neglected terms of order $\hbar^2$ in Equation (\ref{eqn:GS}).

\begin{figure}
  \begin{center}
  \PSbox{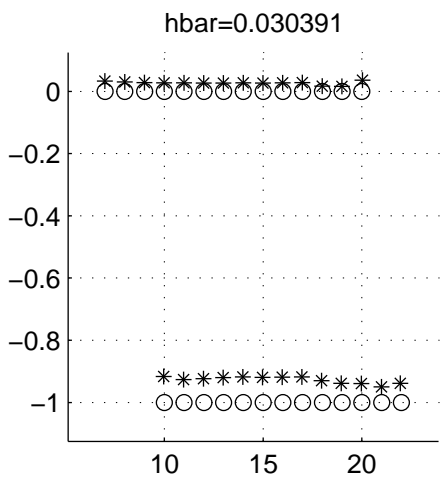}{2.5in}{2.5in}
  \end{center}
  \caption{Lattice points for $\hbar=0.030391$.}
  \label{fig:latt1}
\end{figure}
\begin{figure}
  \begin{center}
  \PSbox{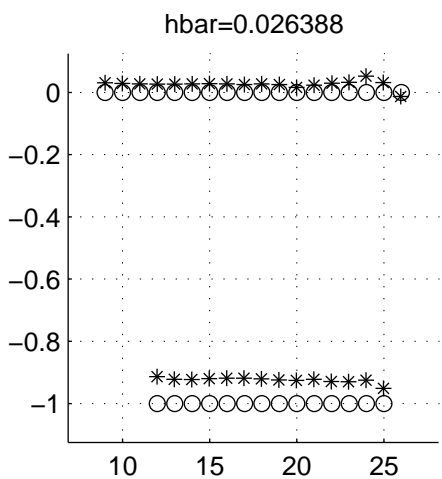}{2.5in}{2.5in}
  \end{center}
  \caption{Lattice points for $\hbar=0.026388$.}
  \label{fig:latt2}
\end{figure}
\begin{figure}
  \begin{center}
  \PSbox{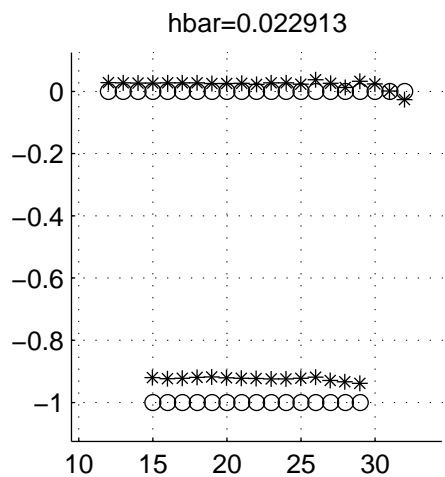}{2.5in}{2.5in}
  \end{center}
  \caption{Lattice points for $\hbar=0.022913$.}
  \label{fig:latt3}
\end{figure}
\begin{figure}
  \begin{center}
  \PSbox{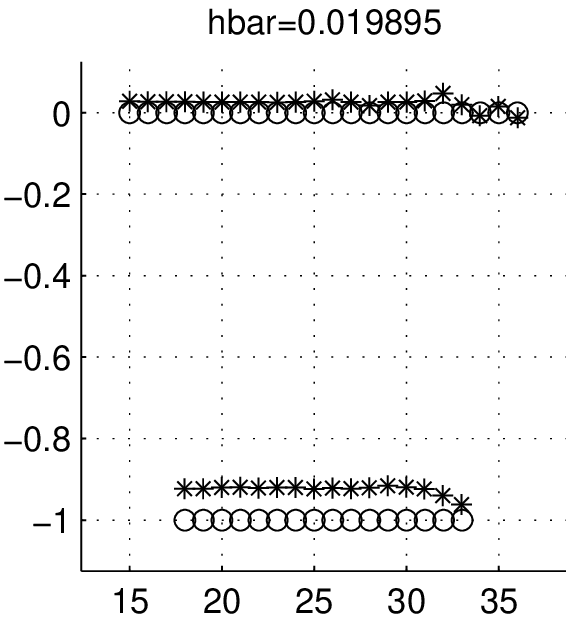}{2.5in}{2.5in}
  \end{center}
  \caption{Lattice points for $\hbar=0.019895$.}
  \label{fig:latt4}
\end{figure}
\begin{figure}
  \begin{center}
  \PSbox{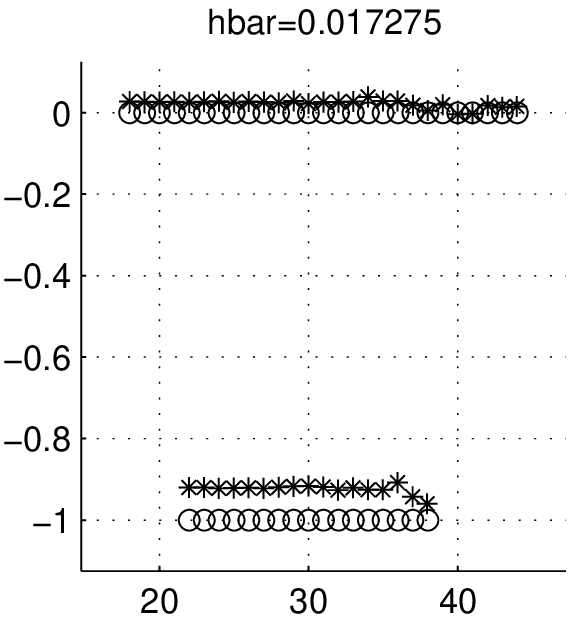}{2.5in}{2.5in}
  \end{center}
  \caption{Lattice points for $\hbar=0.017275$.}
  \label{fig:latt5}
\end{figure}

\section{Conclusions}
Using standard numerical techniqes, one can compute a sufficiently large
number of resonances for the triple gaussian system to verify their
asymptotic distribution in the semiclassical limit
$\hbar\rightarrow{0}$.  This, combined with effective estimates of the
fractal dimension of the classical trapped set, gives strong evidence
that the number of resonances $N_{res}$ in a box
$\interval{E_0,E_1}-i\interval{0,\hbar}$, for sufficiently small
$\abs{E_1-E_0}$ and $\hbar$,
\begin{equation}
  N_{res}\sim\hbar^{-\frac{D\of{K_E}+1}{2}},
\end{equation}
as one can see in Figure {\ref{fig:compare}} and Table
{\ref{tab:compare}}.  Furthermore, the same techniques, when applied to
double gaussian scattering, produce results which agree with rigorous
semiclassical results.  This supports the correctness of our algorithms
and the validity of semiclassical arguments for the range of $\hbar$
explored in the triple gaussian model.  The computation also hints at
more detailed structures in the distribution of resonances: In Figures
{\ref{fig:resplot0}} - {\ref{fig:resplot4}}, one can clearly see gaps
and strips in the distribution of resonances.  A complete understanding
of this structure requires further investigation.

While we do not have rigorous error bounds for the dimension estimates,
the numerical results are convincing.  It seems, then, that the primary
cause for our failure to observe the conjecture in a ``clean'' way is
partly due to the size of $\hbar$: If one could study resonances at much
smaller values of $\hbar$, the asymptotics may become more clear.

\section{Acknowledgments}
Thanks to J.\ Demmel and B.\ Parlett for crucial help with matrix
computations, and to X.\ S.\ Li and C.\ Yang for ARPACK help.  Thanks
are also due to R.\ Littlejohn and M.\ Cargo for their help with bases
and matrix elements, and to F.\ Bonetto for suggesting a practical
method for computing fractal dimensions.  Many thanks to Z.\ Bai, W.\
H.\ Miller, and J. Harrison for helpful conversations, and to the
Mathematics Department at Lawrence Berkeley National Laboratory for
computational resources.  Finally, the author owes much to M.\ Zworski
for inspiring most of this work.

KL was supported by the Fannie and John Hertz Foundation.

\end{document}